\documentclass{article}
\usepackage[paperwidth=7in, paperheight=10in, margin=.875in]{geometry}
\usepackage[backref,colorlinks,linkcolor=red,anchorcolor=green,citecolor=blue]{hyperref}
\usepackage{amsfonts,amssymb}
\usepackage{amsmath}
\usepackage{graphicx}
\usepackage{cite}
\usepackage{enumerate}

\usepackage{comment}

  \newtheorem{theorem}{Theorem}
  \newtheorem{lemma}[theorem]{Lemma}
  \newtheorem{corollary}[theorem]{Corollary}
  \newtheorem{proposition}[theorem]{Proposition}
  
\def\begproof{\noindent{\bf Proof: }}
\def\endproof{\par\rightline{\vrule height5pt width5pt depth0pt}\medskip}

\def\div{\nabla\cdot}
\def\d{\,\mathrm{d}}
\def\eps{\varepsilon}
\def\G{\mathbb{G}}
\def\P{\hbox{\rlap{I}\kern.16em P}}
\def\Q{\hbox{\rlap{\kern.24em\raise.1ex\hbox
      {\vrule height1.3ex width.9pt}}Q}}
\def\M{\hbox{\rlap{I}\kern.16em\rlap{I}M}}
\def\Z{\hbox{\rlap{Z}\kern.20em Z}}
\def\({\begin{eqnarray}}
\def\){\end{eqnarray}}
\def\[{\begin{eqnarray*}}
\def\]{\end{eqnarray*}}

\def\part#1#2{\frac{\partial #1}{\partial #2}}

\def\grad{\nabla}
\def\Norm#1{\left\| #1 \right\|}
\def\bar{\overline}

\def\d{\,\mathrm{d}}
\def\N{\mathbb{N}}
\def\R{\mathbb{R}}

\def\epsilon{\varepsilon}
\def\dx{\,\mathrm{d}x}
\def\ds{\,\mathrm{d}s}

\def\wto{\rightharpoonup}

\def\I{\mathbb{I}}
\def\cc{c}
\def\P{\mathbb{P}}
\def\En{E}
\def\bEn{\bar\En}

\def\bigO{\mathcal{O}}

\newcommand{\Cspace}{C}
\newcommand{\Pspace}{U}

\newcommand{\ptest}{\varphi}
\newcommand{\ptrial}{\varphi}
\newcommand{\Ctest}{\psi}

\newcommand{\ctrial}{\psi}
\newcommand{\dd}[2]{\frac{\partial#1}{\partial #2}}

\newcommand\Vset{\mathbb{V}}
\newcommand\Eset{\mathbb{E}}
\newcommand\E{\mathcal{E}}
\newcommand\CC{\mathcal{C}}
\newcommand\NN{\mathcal{N}}
\newcommand\Qh{\mathbb{Q}^n}
\newcommand\Zh{\mathbb{Z}^n}
\newcommand\Dh{\mathbb{D}^n}
\newcommand\Xh{{X}^n}

\newcommand\Tn{\mathcal{T}^n}

\newcommand{\sz}[1]{{\color{red}#1}}

   \allowdisplaybreaks
\begin{document}
\title{Gradient Flows for the $p$-Laplacian Arising from Biological Network Models: A Novel Dynamical Relaxation Approach}

          \author{
          Jan Haskovec\thanks{Mathematical and Computer Sciences and Engineering Division,
		King Abdullah University of Science and Technology,
		Thuwal 23955-6900, Kingdom of Saudi Arabia, (jan.haskovec@kaust.edu.sa).}
	  \and
          Peter Markowich\thanks{Mathematical and Computer Sciences and Engineering Division,
		King Abdullah University of Science and Technology,
		Thuwal 23955-6900, Kingdom of Saudi Arabia, (peter.markowich@kaust.edu.sa),
		and Faculty of Mathematics, University of Vienna, Oskar-Morgenstern-Platz 1, 1090 Vienna,
		(peter.markowich@univie.ac.at)}
          \and
          Stefano Zampini\thanks{Mathematical and Computer Sciences and Engineering Division,
		King Abdullah University of Science and Technology,
		Thuwal 23955-6900, Kingdom of Saudi Arabia, (stefano.zampini@kaust.edu.sa)}}
        
\date{}

\maketitle

\begin{abstract}
We investigate a scalar partial differential equation model for the formation of biological transportation networks.
Starting from a discrete graph-based formulation on equilateral triangulations, we rigorously derive the corresponding continuum energy functional as the $\Gamma$-limit under network refinement and establish the existence of global minimizers.
The model possesses a gradient-flow structure whose steady states coincide with solutions of the $p$-Laplacian equation.
Building on this connection, we implement finite element discretizations and propose a novel dynamical relaxation scheme that achieves optimal convergence rates in manufactured tests and exhibits mesh-independent performance, with the number of time steps, nonlinear iterations, and linear solves remaining stable under uniform mesh refinement.
Numerical experiments confirm both the ability of the scalar model to reproduce biologically relevant network patterns and its effectiveness as a computationally efficient relaxation strategy for solving $p$-Laplacian equations for large exponents $p$.
\end{abstract}

\vspace{3mm}

\vspace{3mm}


\section{Introduction}
In this paper, we study the scalar partial differential equation system modeling the formation of biological transportation networks,
\(
   - \div [ (\cc + r) \grad u ] &=& S, \label{eq:Poisson} \\
   \part{ \cc}{t} - |\grad u|^2 + \nu \cc^{\gamma-1} &=& 0,  \label{eq:A}
\)
for the scalar pressure field $u=u(t,x)\in\R$ of the fluid transported within the network and the scalar conductance (permeability) $\cc=\cc(t,x)\in\R$,
with metabolic constant $\nu>0$ and metabolic exponent $\gamma>0$.
 We assume the validity of Darcy's law for slow flow in porous media, connecting the flux $q$ of the fluid with the gradient of its pressure via the action of the
 permeability $\P[\cc]$,
\[
   q = - \P[\cc] \grad u.
\]
The total permeability is assumed to be of the form $\P[\cc] := \cc + r$, where
the scalar function $r(x) \geq 0$ describes the isotropic background permeability of the medium.
The source term $S(x)$ in the mass conservation equation \eqref{eq:Poisson}
is to be supplemented as an input datum and is assumed to be independent of time.
We pose \eqref{eq:Poisson}, \eqref{eq:A} on a bounded domain $\Omega\subset\R^d$, $d\geq 1$, with smooth boundary $\partial\Omega$.
We shall consider mixed homogeneous Dirichlet-Neumann boundary conditions for \eqref{eq:Poisson}.
For \eqref{eq:A} we prescribe the initial condition
\( \label{IC_0}
   \cc(t=0,x)=\cc_0(x)\qquad\mbox{for } x\in\Omega,
\)
with $\cc_0(x) \geq 0$ almost everywhere in $\Omega$.

A fundamental structural property of the system \eqref{eq:Poisson}--\eqref{eq:A}
is that it represents the constrained $L^2$-gradient flow associated with the energy functional
\begin{equation} \label{energy}
   \E[\cc] := \int_{\Omega } (\cc+r) |\grad u[\cc]|^2  + \frac\nu\gamma \cc^\gamma \d x,
\end{equation}
where $u[\cc]$ is a solution of the Poisson equation
\eqref{eq:Poisson} subject to boundary conditions.
The first term in \eqref{energy} describes the kinetic (pumping) energy necessary for transporting the medium, while the second term is the metabolic expenditure of the network.
It is well known that for $\gamma\geq 1$ the energy functional is convex, see \cite{HMP22, HV24}.
Based on this fact, we construct a unique solution of  \eqref{eq:Poisson}--\eqref{eq:A}
as the corresponding gradient flow.

A variant of the system \eqref{eq:Poisson}--\eqref{eq:A} was proposed and studied in \cite{Hu-Cai, Hu-Cai-19},
where the local permeability of the porous medium is described in terms
of a vector field.
Detailed mathematical analysis
was carried out in a series of papers \cite{HMP15, HMPS16, bookchapter, AAFM}
and in \cite{Li, Xu2, Xu}, while its various other aspects were studied in \cite{BHMR, HKM2, HMP19, Hu, HV24}.
Another variant with a diagonal permeability tensor was derived formally in \cite{HKM} from the discrete network formation model \cite{Hu}.
The corresponding continuum energy functional was obtained as a limit of a sequence of discrete problems under refinement of the network, with geometry restricted to rectangular networks in two dimensions.
In \cite{HKM2}, the process was carried out rigorously,
in terms of the $\Gamma$-limit of a sequence of discrete energy functionals with respect to the strong $L^2$-topology.
In \cite{HMP22}, this was generalized to general symmetric, positive semidefinite conductivity tensors. The continuum energy functional was obtained as a refinement limit of discrete problems on a sequence of triangular meshes. However, the derivation was only formal.
Finally, in \cite{AH25}, a graphon limit has been derived rigorously, where \eqref{eq:Poisson} is replaced by an integral equation.

A distinctive feature of the scalar network formation system is that its steady states, for metabolic exponents $\gamma>1$, satisfy the $p$-Laplacian equation
\(  \label{eq:p-Laplacian}
|\nabla u|^2 = c^{\gamma -1} \Rightarrow -\nabla  \cdot (|\nabla u|^{p-2} \nabla u) = S, \quad p = \frac{2\gamma}{\gamma - 1} > 2,
\)
when $r=0$, and $\nu=1$. Among the most successful solvers for the $p$-Laplacian are preconditioned gradient descent algorithms \cite{HuangLiLiu2007}, barrier methods \cite{Loisel2020}, and quasi-Newton methods \cite{CaliariZuccher2017}. Within this landscape, relaxed Ka\v{c}anov iterations based on duality arguments \cite{Diening2020, BalciDiening2023} can be considered the gold standard, combining robustness with provable convergence. Semi-implicit formulations arising from gradient flows associated with the $p$-Laplacian energy have also been investigated \cite{BartelsNochetto2018}.
In contrast, the gradient flow proposed here serves as a dynamical relaxation: time integration regularizes the nonlinearity and drives the solution to the $p$-Laplacian steady state, offering an alternative to classical stationary iterations and complementing existing relaxation techniques.

Numerical experiments validate our approach as a competitive alternative to existing relaxation techniques, confirming optimal convergence rates for manufactured solutions and demonstrating robustness for large values of $p$
without the need for adaptive mesh refinement. Numerical evidence indicates that the algorithm exhibits mesh-independent convergence, in the sense that the number of required time steps, nonlinear iterations, and inner linear solver steps remains essentially unaffected by mesh resolution.

This paper is structured as follows. In Section \ref{sec:derivation} we rigorously derive the energy functional \eqref{energy} with $\gamma>1$ from the discrete model \cite{Hu-Cai-19, Hu-Cai} in the refinement $\Gamma$-limit of a sequence of triangulations. We also construct its global minimizers in the class of uniformly Lipschitz continuous functions as limits of minimizers of the discrete problem.
In Section \ref{sec:numerics} we numerically investigate the system \eqref{eq:Poisson}--\eqref{eq:A}. First, in Section \ref{sec:sequence_formation} we consider the case $\gamma\in (0,1)$ and demonstrate its ability to generate two-dimensional network patterns. Then, in Section \ref{subsec:pLaplacian} we solve \eqref{eq:Poisson}--\eqref{eq:A} with $\gamma>1$ until equilibrium to obtain solutions of the $p$-Laplacian equation \eqref{eq:p-Laplacian}. 


\section{Derivation from the discrete model in two dimensions}\label{sec:derivation}

The goal of this Section is to provide a derivation of the continuum energy functional \eqref{energy} with $\gamma>1$ as a limit of discrete energy functionals \cite{Hu-Cai-19, Hu-Cai} posed on a sequence of equilateral
triangulations of a bounded polygonal domain $\Omega\subset \R^2$. For simplicity and without loss of generality, we will assume that $r>0$ is a constant, and we will consider homogeneous Neumann boundary conditions.

We construct a sequence of undirected graphs $\G^n=(\Vset^n,\Eset^n)$, $n\in\N$,
consisting of finite sets of vertices $\Vset^n$ and edges $\Eset^n$,
which are all equilateral triangulations of $\Omega$ with edge length $h>0$ approaching zero as $n\to\infty$.
For any $n\in\N$ we shall denote $\Tn$ the set of all the triangles in the corresponding equilateral triangulation.
One possible way is to iteratively refine each triangle into four equilateral triangles with half-edge lengths.
This would imply $h = \bigO(2^{-n})$ and the number of triangles (and edges, nodes) proportional to $2^{2n}$.
However, the methods presented in this paper
work for any sequence of equilateral triangulations such that $h\to 0$ as $n\to\infty$. In fact, our approach can be generalized to quasi-uniform triangulations. However, this would severely increase the technicality of the exposition. We therefore choose to restrict to the equilateral case.

Each edge $(i,j)\in\Eset^n$, connecting the vertices $i$ and $j\in\Vset^n$, has an associated
conductivity denoted by $\CC_{ij}=\CC_{ji}\geq 0$.
Moreover, each vertex $i\in\Vset^n$ has the pressure $U_i\in\R$ of the material flowing through it.
The local mass conservation in each vertex is expressed in terms of the Kirchhoff law
\begin{align}\label{eq:kirchhoff}
   -\sum_{j\in \NN(i)} (\CC_{ij}+r)\frac{U_j-U_i}{h}=S_i\qquad \text{for all~}i\in \Vset^n,
\end{align}
where $\NN(i)$ denotes the set of vertices connected to $i\in\Vset^n$ through an edge.
Moreover, $S=(S_i)_{i\in\Vset^n}$ is the prescribed vector of strengths of flow sources ($S_i>0$) or sinks ($S_i<0$).
Given the vector of conductivities $\CC=(\CC_{ij})_{(i,j)\in\Eset^n}$,  the Kirchhoff law \eqref{eq:kirchhoff}
is a linear system of equations for the vector of pressures $U=(U_i)_{i\in\Vset^n}$.
A necessary condition for its solvability is the global mass balance
\(   \label{mass_balance}
   \sum_{i\in\Vset^n} S_i=0,
\)
whose validity we assume in the sequel.
Note that the solution $U=(U_i)_{i\in\Vset^n}$ is unique up to an additive constant.

The conductivities $\CC_{ij}$ are subject to a minimization process with the energy functional
\begin{align}\label{eq:energydisc}
   {\En}^n[\CC] := \, h \sum_{(i,j)\in\Eset^n} (\CC_{ij}+r) \frac{(U_j-U_i)^2}{h^2} + \frac\nu\gamma \CC_{ij}^\gamma ,
\end{align}
where $\CC = (\CC_{ij})_{(i,j)\in\Eset^n}$ denotes the vector of edge conductivities.
The first term in the summation describes the (physical) pumping power necessary for transporting the material through the network,
while the second term describes the metabolic cost of maintaining the network structure.

For convenience, we shall use both notations $(i,j)\in\Eset^n$ to address edges by the indices
of their adjacent vertices, and $e_{ij}\in\Eset^n$ to address the linear segment in $\R^2$ connecting
them (i.e., the set of all convex combinations of the vertex coordinates $x_i$, $x_j \in\R^2$).
For each edge $e_{ij} \in\Eset^n$ we denote by $\diamond_{ij}^n\subset\R^2$
the adjacent closed diamond, i.e.,
\[
   \diamond_{ij}^n := \left\{ x\in\Omega;\, \mbox{dist}(x,e_{ij}) \leq \mbox{dist}(x,\tilde e) \mbox{ for all } \tilde e\in\Eset^n \right\}.
\]
Since the graph $\G^n=(\Vset^n,\Eset^n)$ represents an equilateral triangulation of $\Omega$, the diamonds are parallelograms with angles $60^{\circ} $ and $120^{\circ} $.
For edges $e_{ij}$ belonging to the boundary $\partial\Omega$ the definition of $\diamond_{ij}$
has to be adjusted accordingly, such that the corresponding diamond is also a parallelogram with angles $60^{\circ} $ and $120^{\circ}$.
For each vertex $i\in\Vset^n$ we denote $\Diamond^n_i$ the union of adjacent diamonds,
\[
   \Diamond^n_i := \bigcup_{j\in \NN(i)} \diamond_{ij}^n,
\]
where $\NN(i)$ again denotes the set of vertices connected to $i$.

The main idea of our construction is to map the vector of discrete edge conductivities 
onto the set of piecewise constant scalar fields on $\Omega$ through the mapping $\Qh$
defined as
\(    \label{Qh}
   \Qh[\CC](x) := \CC_{ij}  \qquad\mbox{for all } x\in \diamond_{ij}^n.
\)
$\Qh[\CC]$ is a constant scalar on each diamond $\diamond_{ij}^n$, taking the value $\CC_{ij}$.
We also introduce the operator $\Zh$
that maps the edge conductivities 
onto piecewise constant scalar fields on the triangles $T\in\Tn$,
\(    \label{Zh}
   \Zh[\CC](x) := \frac13 \sum_{(i,j)\in T} \CC_{ij}  \qquad\mbox{for all } x\in T,
\)
where the sum goes through the three edges constituting the triangle $T$.
We also introduce the piecewise constant vector field $\Xh\in\R^2$, defined as
\(    \label{Xh}
   \Xh(x) :=  \frac{x_i-x_j}{|x_i-x_j|} \qquad\mbox{for all } x\in \diamond_{ij}^n,
\)
where we recall that $x_i\in\Omega$ and, resp., $x_j\in\Omega$ denote the spatial coordinates of the vertices $i\in\Vset^n$ and, resp., $j\in\Vset^n$.

Our strategy is, by means of the mappings $\Qh$ and $\Zh$, to establish a connection between a properly rescaled
version of the discrete energy functional \eqref{eq:energydisc} and its continuum counterpart \eqref{energy}.
Following the strategy of \cite{HKM2}, the connection will be established through an "intermediate"
functional $\E^n$, where the pressure is a solution of a finite element discretization of the Poisson
equation on the triangulation $\Tn$. We shall call this functional the \emph{semi-discrete energy functional}. As the first step we study the connection between the Kirchhoff law \eqref{eq:kirchhoff} and the Poisson equation
\eqref{eq:Poisson}.

\subsection{The Kirchhoff law and the Poisson equation}\label{subsec:Poisson}

We consider the weak formulation of the Poisson equation \eqref{eq:Poisson} with the nonnegative permeability $\cc\in L^\infty(\Omega)$ and homogeneous Neumann boundary condition,
\(  \label{Poisson:weak}
   \int_\Omega (\cc + r) \grad u\cdot\grad\psi \,\d x =
   \int_\Omega S\psi\,\d x
   \qquad
   \mbox{for all } \psi\in L^2(\Omega).
\)
Throughout Section \ref{sec:derivation} we impose the global mass balance
\(  \label{ass:S}
   \int_\Omega S(x) \d x = 0.
\)
We first establish the solvability of \eqref{Poisson:weak} with a square-integrable permeability kernel $\cc\in L^2(\Omega)$. In the sequel we denote $H^1_0(\Omega)$ the set of functions in the Sobolev space $H^1(\Omega)$ with zero mean.

\begin{proposition}\label{lem:Poisson}
Let $r>0$, $S \in L^2(\Omega)$ satisfying the global mass balance \eqref{ass:S}, and $\cc \in L^2(\Omega)$ with $\cc\geq 0$ almost everywhere on $\Omega$.
Then  there exists a unique $u \in H^1_0(\Omega)$ verifying \eqref{Poisson:weak}
for all test functions $\psi\in L^\infty(\Omega)$.
Moreover, we have
\(   \label{est:p}
     \Norm{\grad u}_{L^2(\Omega)}  \leq  \frac{C_\Omega}{r} \Norm{S}_{L^2(\Omega)},
\)
and
\(  \label{eq:testPoisson}
   \int_\Omega (\cc+r) |\grad u|^2 \d x = \int_\Omega S u \,\d x \leq \frac{C_\Omega^2}{r} \Norm{S}^2_{L^2(\Omega)},
\)
where $C_\Omega$ is the Poincar\'{e} constant of the domain $\Omega$.
\end{proposition}

The proof is obtained as a slight modification of the proof of \cite[Lemma 6]{HMP15} or \cite[Lemma 5.5]{AH25}.
Now, for the nonnegative permeability $\cc\in L^\infty(\Omega)$
we consider the modified Poisson equation
\(  \label{eq:PFEMstrong}
   - 2\grad\cdot \bigl( (\cc+r) (\Xh\otimes\Xh)\grad u \bigr) = S,
\)
with $\Xh$ defined in \eqref{Xh} and the symbol $\otimes$ denoting the standard tensor product. The factor $2$ is introduced because, as we shall see, $\Xh\otimes\Xh\wto\frac12\mathbb{I}$ as $n\to\infty$.
In the sequel we shall work with the first-order $H^1$ finite element approximation of \eqref{eq:PFEMstrong},
\(   \label{eq:PFEM}
   2\int_\Omega (\cc+r)\, \grad \psi \cdot (\Xh\otimes\Xh) \grad u\, \d x = \int_\Omega S \psi \,\d x
      \qquad \mbox{for  all } \psi \in Y_1^n.
\)
Here $Y_1^n$ denotes the space of continuous, piecewise linear functions on $\Tn$,
\[
   Y_1^n := \left\{ u\in C^0(\Omega);\; u \mbox{ linear on each } T\in \Tn \right\}.
\]
Moreover, we denote $Y_0^n$ the space of bounded, piecewise constant functions on the triangulation $\Tn$,
\[
   Y_0^n := \left\{ u\in L^\infty(\Omega);\; u \mbox{ constant on each } T\in \Tn \right\}.
\]
Existence of solutions of \eqref{eq:PFEM} is established by a standard application of the Lax-Milgram theorem, where coercivity of the respective bilinear form with $r>0$ follows from formula \eqref{eq:T} below.
The following proposition establishes a connection between the solution of \eqref{eq:PFEM} with $\cc:=\Qh[\CC]$ and the Kirchhoff law \eqref{eq:kirchhoff}.

\begin{proposition}\label{prop:FEM}
Let $S\in L^2(\Omega)$ verify the global mass balance assumption \eqref{ass:S} and fix $n\in\N$.
Let $\CC_{ij}\geq 0$ for all $(i,j)\in\Eset^n$ and let $u^n:=u[\Qh[\CC]]\in Y_1^n$ be the solution of the discrete Poisson equation \eqref{eq:PFEM} with conductivity $\cc:=\Qh[\CC]$. For $i\in\Vset^n$ denote $U_i^n := u^n(x_i)$.

Then we have, for all $i\in\Vset^n$,
\(  \label{Kresc}
     - \sum_{j\in \NN(i)} (\CC_{ij} + r) \frac{U_j^n-U_i^n}{h} = S_i^n,
\)
where we denoted
\(  \label{Si}
   S^n_i := \frac{h}{2\mbox{vol}\; (\diamond_h)}\int_\Omega S \psi^n_i \,\d x \qquad\mbox{for all } i\in\Vset^n,
\)
and $\psi_i^n\in Y_1^n$ is a piecewise linear function supported on $\Diamond_i$, taking the vertex values
$\psi_i^n(x_i)=1$ and $\psi_i^n(x_j)=0$ for all $j\neq i$.
\end{proposition}

The proof is a minor modification of the proof of \cite[Proposition 2.1]{HMP22}.
Note that \eqref{Kresc} is the Kirchhoff law \eqref{eq:kirchhoff} with right-hand side given by \eqref{Si}. 
We have
\[
   \frac{2\mbox{vol}(\diamond_h)}{h}\sum_{i\in\Vset^n} S^n_i &=& \sum_{i\in\Vset^n} \int_\Omega S(x) \psi^n_i(x) \d x  \\
      &=& \sum_{T\in\Tn} \int_T S(x) \sum_{i\in\Vset^n} \psi^n_i(x) \d x  \\
      &=& \int_\Omega S(x) \d x,
\]
where we used the fact that the sum of barycentric coordinates on each triangle is constant unity,
i.e., $\sum_{i\in\Vset^n} \psi^n_i(x) \equiv 1$ for all $x\in\Omega$.
Consequently, assumption \eqref{ass:S} implies the discrete mass balance \eqref{mass_balance}.

Moreover, for equilateral triangulations we have
\[
   \frac{h}{2\mbox{vol}(\diamond_h)}= \frac{\sqrt{3}}{h},
   \qquad
   \int_\Omega (\psi_i^n)^2 \d x = \frac{\sqrt{3}}{24} h^2.
\]
Consequently, by the Cauchy-Schwarz inequality,
\[
   (S_i^n)^2 \leq \frac{\sqrt{3}}{8} \int_{\Diamond_i} S(x)^2 \d x,
\]
so that
\(\label{est:Si}
   \sum_{i\in\Vset^n} (S_i^n)^2 \leq \frac{3\sqrt{3}}{8} \int_\Omega S(x)^2 \d x.
\)

To facilitate the limit passage $n\to\infty$, we shall work with the discrete Poisson equation \eqref{eq:PFEM}
with $\cc:=\Zh[\CC]$. 
We have the following result.

\begin{lemma} \label{lem:T}
Fix $n\in\N$.
For any $u,v \in Y_1^n(\Omega)$ and any equilateral triangle $T\in\Tn$ we have
\(  \label{eq:T}
   \int_T \grad u\cdot (\Xh\otimes\Xh) \grad v \,\d x = \frac12 \int_T \grad u\cdot\grad v \,\d x.
\)
\end{lemma}

\begproof
Let us denote the three edges of the triangle $T$ by $e_1$, $e_2$, $e_3$
and, without loss of generality, assume $|e_1| = |e_2| = |e_3| = 1$.
Since both $\grad u$ and $\grad v$ are, by definition, constant on the triangle, with \eqref{Xh} we have
\[
   \int_T \grad u\cdot (\Xh\otimes\Xh) \grad v \,\d x &=&
      \grad u\cdot \left( \int_T \Xh\otimes\Xh \,\d x \right) \grad v  \\
      &=&
      \grad u\cdot \left( \frac{1}{4\sqrt 3} \sum_{i=1}^3 e_i\otimes e_i \right) \grad v.
\]
By rotational symmetry, the matrix $\sum_{i=1}^3 e_i\otimes e_i$ is a scalar multiple of the identity.
Since
\[
   \mbox{trace} \left( \sum_{i=1}^3 e_i\otimes e_i \right) = \sum_{i=1}^3 |e_i|^2 = 3,
\]
we conclude $\sum_{i=1}^3 e_i\otimes e_i = \frac{3}{2} \I$.
Consequently,
\[
   \int_T \grad u\cdot (\Xh\otimes\Xh) \grad v \,\d x =
      \frac{\sqrt{3}}{8} \grad u\cdot \grad v =
      \frac12 \int_T \grad u\cdot\grad v \,\d x.
\]
\endproof

Lemma \ref{lem:T} allows us to write the left-hand side of \eqref{eq:PFEM} with $\cc:=\Zh[\CC]$ as
\[
   \int_\Omega (\Zh[\CC]+r)\, \grad \psi \cdot (\Xh\otimes\Xh) \grad u\, \d x
   &=& \sum_{T\in\Tn}  (\Zh[\CC]_T+r) \int_T \grad \psi \cdot (\Xh\otimes\Xh)  \grad u \,\d x  \\
   &=& \frac12 \int_\Omega (\Zh[\CC]+r) \grad \psi\cdot  \grad u \, \d x,
\]
where $\Zh[\CC]_T$ denotes the constant value of $\Zh[\CC]$ on the triangle $T\in\Tn$.
Consequently, \eqref{eq:PFEM} with $\cc:=\Zh[\CC]$ can be equivalently formulated as
\(  \label{PZweak}
      \int_\Omega ( \Zh[\CC] + r ) \grad \psi \cdot \grad u \,\d x = \int_\Omega S \psi \, \d x
      \qquad \mbox{for  all } \psi \in Y_1^n.
\)
This formulation facilitates the passage to the $\Gamma$-limit as $n\to+\infty$
since it avoids the oscillatory term $\Xh\otimes\Xh$ that converges only weakly to $\frac12\I$.

To estimate the difference between the scalar fields $\Qh[\CC]$ and $\Zh[\CC]$, let us for any $T\in\Tn$ denote
\( \label{def:Dh}
   \Dh[\CC](T):= \max\{ |\CC_{ij}-\CC_{ik}|,|\CC_{ik}-\CC_{jk}|,|\CC_{ij}-\CC_{jk}|\}
\)
where $(i,j,k)$ are the vertices constituting $T$.

\begin{lemma}\label{lem:23CC}
For any $\CC\in\R^{|\Eset^n|}$ we have
\(\label{eq:23CC}
   \Norm{\Qh[\CC]-\Zh[\CC]}_{L^\infty(\Omega)}
   \leq \frac23 \max_{T\in\Tn} \Dh[\CC](T).
\)
\end{lemma}

\begproof
Let us fix any triangle $T\in\Tn$ and denote its
vertices $i$, $j$, $k\in\Vset^n$.
For the half-diamond $\diamond_{ij}\cap T$, we have
\[
   \Qh[\CC] \equiv \CC_{ij},\qquad
   \Zh[\CC] \equiv \frac13 (\CC_{ij}+\CC_{jk}+\CC_{ik}).
\]
Consequently,
\[
   \left| \Qh[\CC^n] - \Zh[\CC^n] \right|
   =
   \frac13 \left| 2\CC_{ij} - \CC_{jk} - \CC_{ik} \right|
   \leq
   \frac23 \max\left\{ \left| \CC_{ij} - \CC_{jk} \right|,
   \left| \CC_{ij} - \CC_{ik} \right| \right\}.
\]
Repeating the estimate for the other two half-diamonds constituting $T$, and taking maximum over all $T\in\Tn$, we obtain \eqref{eq:23CC}.
\endproof

\subsection{The semi-discrete energy functionals}\label{subsec:semi-disc}

For the nonnegative permeability $\cc\in L^\infty(\Omega)$
we introduce the semi-discrete energy functional
\(    \label{def:EEh}
    \E^n[\cc] = \int_\Omega 2\left(\cc + r\right) \left| X^n\cdot \grad u[\cc] \right|^2 + \frac\nu\gamma \cc^\gamma \d x,
\)
where $u[\cc]\in Y_1^n$ the unique weak solution of the discrete Poisson equation \eqref{eq:PFEM}.

Moreover, we introduce the rescaled version of the discrete energy functional \eqref{eq:energydisc},
\(  \label{def:en:disc:resc}
   \bar \En^n[\CC] =  {\mbox{vol}(\diamond_h)} \sum_{(i,j)\in\Eset^n}  2\left(\CC_{ij}+ r \right) \frac{(U_i - U_j)^2}{h^2} + \frac\nu\gamma \CC_{ij}^\gamma ,
\)
where $U=U[\CC]$ is a solution of the rescaled Kirchhoff law \eqref{Kresc}.
The following proposition, which is a slight modification of \cite[Proposition 2.2]{HMP22},
establishes the equality of $\E^n$ evaluated at the piecewise constant function $\Qh[\CC]$
to the value of $\bar E^n[\CC]$.

\begin{proposition}\label{prop:En}
For any $n\in\N$ and any vector of nonnegative conductivities $\CC = (\CC_{ij})_{(i,j)\in\Eset^n}$ we have
\(
   \bEn^n[\CC] = \E^n\left[\Qh[\CC] \right],
\)
with the rescaled discrete energy functional $\bEn^n$ defined in \eqref{def:en:disc:resc} and
the semi-discrete energy $\E^n$ given by \eqref{def:EEh}.
\end{proposition}

We now estimate the difference of energies calculated with $\cc:=\Qh[\CC]$ and $\cc:=\Zh[\CC]$.

\begin{lemma}\label{lem:EstDiffE}
Fix $n\in\N$.
For any $\CC\in\R^{|\Eset^n|}$ we have
\(\label{est:DiffE}
   &&\left| \E^n[\Qh[\CC]] - \E^n[\Zh[\CC]] \right| \\
   &&\qquad\leq
   \frac23 \max_{T\in\Tn} \Dh[\CC](T)
   \left( 2
   \Norm{\grad u^\mathbb{Q}}_{L^2(\Omega)}
   \Norm{\grad u^\mathbb{Z}}_{L^2(\Omega)}
   + |\Omega|\nu \Norm{\CC}_{\ell_\infty}^{\gamma-1}
   \right),
   \nonumber
\)
where we denoted $u^\mathbb{Q}\in Y_1^n$ the solution
of the discrete Poisson equation \eqref{eq:PFEM}
with permeability $\Qh[\CC]$,
and $u^\mathbb{Z}\in Y_1^n$ the solution of \eqref{eq:PFEM} with permeability $\Zh[\CC]$.
\end{lemma}

\begproof
Using $u^\mathbb{Q}$ as a test function in the discrete Poisson equation \eqref{eq:PFEM} for $u^\mathbb{Q}$ yields
for the kinetic part of \eqref{def:EEh},
\[
   \E^n_{\mathrm{kin}}[\Qh[\CC]] =
   2\int_\Omega \Qh[\CC] |X^n\cdot\grad u^\mathbb{Q}|^2 \d x
   = \int S u^\mathbb{Q} \d x.
\]
Analogously, we have
\[
   \E^n_{\mathrm{kin}}[\Zh[\CC]] =
   2\int_\Omega \Zh[\CC] |X^n\cdot\grad u^\mathbb{Z}|^2 \d x
   = \int S u^\mathbb{Z} \d x.
\]
Therefore,
\[
   \E^n_{\mathrm{kin}}[\Qh[\CC]] - \E^n_{\mathrm{kin}}[\Zh[\CC]]
   = \int_\Omega S \left(u^\mathbb{Q} - u^\mathbb{Z}\right) \d x.
\]
Now, using $u^\mathbb{Z}$ as a test function in the discrete Poisson equation for $u^\mathbb{Q}$ and vice versa, and subtracting the resulting identities, we obtain
\[
   2\int_\Omega \left( \Qh[\CC] - \Zh[\CC] \right)
   (\Xh\cdot \grad u^\mathbb{Q}) (\Xh\cdot \grad u^\mathbb{Z}) \d x
   = \int_\Omega S \left(u^\mathbb{Q} - u^\mathbb{Z}\right) \d x.
\]
We combine the last two identities
and use the Cauchy-Schwarz inequality to estimate
\[
   && \left| \int_\Omega \left( \Qh[\CC] - \Zh[\CC] \right)
   (\Xh\cdot \grad u^\mathbb{Q}) (\Xh\cdot \grad u^\mathbb{Z}) \d x
   \right| \\
   &&\qquad \leq
   \Norm{\Qh[\CC]-\Zh[\CC]}_{L^\infty(\Omega)}
   \Norm{\grad u^\mathbb{Q}}_{L^2(\Omega)}
   \Norm{\grad u^\mathbb{Z}}_{L^2(\Omega)},
\]
where we also used the uniform bound $|\Xh|\leq 1$.

To estimate the difference of the metabolic parts of \eqref{def:EEh}, we use the inequality $|a^\gamma - b^\gamma| \leq \gamma\max\{a,b\}^{\gamma-1} |a-b|$, which holds for $a$, $b\geq 0$ and $\gamma\geq 1$. We obtain
\[
   \frac\nu\gamma \int_\Omega \left|\Qh[\CC]^\gamma - \Zh[\CC]^\gamma \right| \d x
   \leq
   |\Omega|\nu \Norm{\CC}_{\ell_\infty}^{\gamma-1} \Norm{\Qh[\CC]-\Zh[\CC]}_{L^\infty(\Omega)}.
\]
We conclude by using \eqref{eq:23CC}.
\endproof

With the uniform $H^1$-bound \eqref{est:p} on the solutions of the Poisson equation \eqref{eq:PFEM} provided by Proposition \ref{lem:Poisson}, we readily obtain the following corollary.

\begin{corollary}\label{lem:epsn}
Let $\CC^n \in \R^{|\Eset^n|}$, $n\in\N$, be a sequence of nonnegative conductivities such that
\(\label{eq:coro}
   \sup_{n\in\N} \Norm{\CC^n}_{\ell_\infty} <\infty
   \qquad\mbox{and}\qquad
   \lim_{n\to\infty} \max_{T\in\Tn} \Dh[\CC^n](T) = 0.
\)
Then
\[
  \lim_{n\to\infty} \left| \E^n[\Qh[\CC^n]] - \E^n[\Zh[\CC^n]] \right| = 0.
\]
\end{corollary}

\subsection{The continuum energy functional}
In this section we study the properties of the continuum energy functional \eqref{energy} and its relation to its semi-discrete version \eqref{def:EEh}.
We start by proving the continuity of \eqref{energy} with respect to the norm topology of $L^\omega(\Omega)$, $\omega:=\max\{2,\gamma\}$.

\begin{lemma}\label{lem:strong}
Let $S \in L^2(\Omega)$ satisfy \eqref{ass:S} and fix $r>0$.
Let the sequence $(\cc^n)_{n\in\N} \subset L^\omega(\Omega)$
converge to $\cc\in L^\omega(\Omega)$ in the norm topology of $L^\omega(\Omega)$ as $n\to \infty$.
Then
\[
   \E[\cc] = \lim_{n\to\infty} \E[\cc^n].
\]
\end{lemma}

\begproof
Let $(u^n)_{n\in\N}\subset L^2_0(\Omega)$ be a sequence of
weak solutions of the Poisson equation \eqref{Poisson:weak} with permeability $\cc^n\in L^\omega(\Omega)$,
provided by Proposition \ref{lem:Poisson}.
Identity \eqref{eq:testPoisson} gives
\(  \label{eq:testPoisson2}
   \int_\Omega (\cc^n + r) |\grad u^n|^2 \d x = \int_\Omega S u^n \d x.
\)
Due to the uniform bound \eqref{est:p},
there exists $u\in H^1_0(\Omega)$, a weak limit of (a subsequence of) $u^n$,
verifying \eqref{Poisson:weak} for all test functions $\varphi\in L_0^\infty(\Omega)$.
Since $\cc^n\to \cc$ strongly in $L^2(\Omega)$, we can pass to the limit in the weak formulation \eqref{Poisson:weak}
which establishes $u$ as its weak solution with permeability $\cc$
and test functions $\varphi\in L^\infty(\Omega)$.
Again, \eqref{eq:testPoisson} gives
\[
   \int_\Omega (\cc+r) |\grad u|^2 \d x = \int_\Omega S u \, \d x.
\]
Combining with \eqref{eq:testPoisson2} we have
\[
   \lim_{n\to\infty} \int_\Omega (\cc^n+r) |\grad u^n|^2 \d x =
   \lim_{n\to\infty}  \int_\Omega S u^n \d x =  \int_\Omega S u \d x =
   \int_\Omega (\cc+r) |\grad u|^2 \d x,
\]
which proves the strong continuity of the kinetic part of the energy functional $\En$.
The metabolic part is a multiple of the $L^\gamma$-norm of $\cc$, consequently, its strong continuity is obvious.
\endproof

A major observation is that for $\gamma\geq 1$
the energy functional $\E$ given by \eqref{energy} is convex. A proof can be found in \cite[Proposition 3.2]{HMP22} or \cite[Lemma 6.4]{AH25}.
Then, since $\E$ is strongly continuous on $L^\omega(\Omega)$ due to Lemma \ref{lem:strong},
an application of the Mazur's lemma \cite{Rockafellar}
establishes its weak lower semicontinuity, i.e.,
\(  \label{eq:wlsc}
   \E[\cc] \leq \liminf_{n\to \infty} \E[\cc^n]
\)
whenever $(\cc^n)_{n\in\N} \subset L^\omega(\Omega)$ converges weakly in $L^\omega(\Omega)$ to $\cc$.

Next, we establish an inequality between the semi-discrete and continuum energy functionals
restricted to piecewise constant permeabilities on the triangles $T\in\Tn$.

\begin{lemma}\label{lem:ineqE}
For any $n\in\N$ and nonnegative $\cc^n\in Y^n_0$ we have
\(   \label{EnnEn}
   \E^n[\cc^n] \leq \E[\cc^n],
\)
with the semi-discrete energy functional $\E^n$ given by \eqref{def:EEh} and $\E$ given by \eqref{energy}.
\end{lemma}

\begproof
Let us recall the weak formulation of the discrete Poisson equation \eqref{eq:PFEM} with test function $\psi\in Y_1^n$,
\[
   2\int_\Omega (\cc^n+r)\, \grad \psi \cdot (\Xh\otimes\Xh) \grad u^n\, \d x = \int_\Omega S \psi \,\d x.
\]
Lemma \ref{lem:T} gives
\[
   2\int_\Omega (\cc^n+r)\, \grad \psi \cdot (\Xh\otimes\Xh) \grad u^n\, \d x
   &=&
     2\sum_{T\in\Tn}  \left(\cc_T^n + r \right)  \int_T (\Xh\cdot \grad u^n)(\Xh\cdot\grad\psi) \,\d x \\
     &=&
       \sum_{T\in\Tn} |T| \left(\cc_T^n + r \right)  (\grad u^n)_T \cdot (\grad\psi)_{T},
\]
where $\cc_T^n$ denotes the constant value of $\cc^n$ on the triangle $T\in\Tn$,
and similarly for $(\grad\psi)_T$ and $(\grad u^n)_T$.

Let $u\in H^1_0(\Omega)$ be the unique weak solution of the Poisson equation \eqref{eq:Poisson} with permeability $c^n$.
Using again $\psi\in Y_1^n$ as a test function gives
\[
   \int_\Omega (\cc^n +r) \grad u \cdot \grad\psi \,\d x =
         \sum_{T\in\Tn} \left(\cc_T^n + r \right) (\grad\psi)_T \cdot \int_{T} \grad u \,\d x.
\]
Due to the uniqueness (in the respective spaces of functions with vanishing mean) of the solutions of \eqref{eq:PFEM} and \eqref{eq:Poisson} we have
\[
   (\grad u^n)_T = \frac{1}{|T|} \int_T \grad u \, \d x \qquad\mbox{for all } T\in\Tn.
\]
Then we have for the kinetic parts of the semi-discrete and continuum energies
\[
    \int_\Omega (\cc^n+r) |\grad u^n|^2 \d x 
      &=& \sum_{T\in\Tn} |T| (\cc^n_T+r) |(\grad u^n)_T|^2 \\
      &=& \sum_{T\in\Tn} |T| (\cc^n_T+r)  \left| \frac{1}{|T|} \int_T \grad u \, \d x \right|^2 \\
      &\leq& \sum_{T\in\Tn} (\cc^n_T+r)  \int_T |\grad u|^2 \d x \\
      &=& \int_\Omega (\cc^n+r) |\grad u|^2 \d x,
\]
where we used the Jensen inequality. Since the metabolic part $\frac{\nu}{\gamma} \int_\Omega (c^n)^\gamma \d x$
is identical for both the semi-discrete and continuum energies, we conclude \eqref{EnnEn}.
\endproof

\begin{lemma} \label{lem:Enn}
Let $(\cc^n)_{n\in\N}$, be a sequence of piecewise constant nonnegative conductivities $\cc^n \in Y_0^n$
such that $sup_{n\in\N} \Norm{c^n}_{L^\infty(\Omega)} <\infty$.
Then
\(
   \lim_{n\to\infty} \left| \E^n[\cc^n] - \E[\cc^n]\right| = 0.
\)
\end{lemma}

\begproof
Using the solution $u^n\in H^1_0(\Omega)$ of \eqref{Poisson:weak} as a test function, we obtain
\[
   \E_{\mathrm{kin}}[\cc^n] = \int_\Omega (\cc^n +r) |\grad u^n|^2 \d x = \int_\Omega S u^n\,\d x.
\]
On the other hand, using the solution $u_h^n\in Y^n_1(\Omega)$ of \eqref{eq:PFEM} as a test function, we obtain
\[
   \E^n_{\mathrm{kin}}[\cc^n] = \int_\Omega 2(\cc^n+r) |X^n\cdot\grad u_h^n|^2 \d x = \int_\Omega S u_h^n\,\d x.
\]
Consequently,
\[
   \left| \E^n[\cc^n] - \E[\cc^n]\right|
   &=&
   \left| \int_\Omega S \left(u^n - u_h^n\right) \, \d x \right| \\
   &\leq&
   \Norm{S}_{L^2(\Omega)} \Norm{u^n - P^n[u^n]}_{L^2(\Omega)} +
   \left|\int_\Omega S \left(P^n[u^n] - u_h^n\right) \,\d x \right|,
\]
where $P^n : H^1(\Omega) \to Y^n_1$ is an $H^1$-stable interpolator, i.e., Scott–Zhang or Clement. 
Note that due to the uniform boundedness of $u^n$ in $H^1_0(\Omega)$ provided by \eqref{est:p}, we have
\( \label{H1stable1}
   \lim_{n\to\infty} \Norm{u^n - P^n[u^n]}_{H^1(\Omega)} = 0.
\)
Moreover, using $P^n[u^n]\in Y^n_1$ as a test function in the discrete Poisson equation \eqref{eq:PFEM} for $u_h^n$,
we have
\[
   \int_\Omega 2(\cc^n +r) (X^n\cdot\grad u_h^n) (X^n\cdot \grad P^n[u^n]) \,\d x
   &=&
   \int_\Omega (\cc^n +r) \grad u_h^n \cdot \grad P^n[u^n] \,\d x \\
   &=&
   \int_\Omega S P^n[u^n] \,\d x,
\]
where we used Lemma \ref{lem:T} in the first line, recalling that $\cc^n \in Y_0^n$.
On the other hand, using $u_h^n\in Y^n_1$ as a test function in the Poisson equation \eqref{Poisson:weak} for $u^n$, we obtain
\[
   \int_\Omega (\cc^n +r) \grad u^n\cdot \grad u^n_h \,\d x
   =
   \int_\Omega S u^n_h \,\d x.
\]
Consequently,
\[
   \left|\int_\Omega S \left(P^n[u^n] - u_h^n\right) \,\d x \right|
   &=&
   \left| \int_\Omega (\cc^n +r) \grad u_h^n\cdot\grad\left( P^n[u^n] - u^n\right)\, \d x \right| \\
   &\leq&
   \sup_{n\in\N} \Norm{\cc^n+r}_{L^\infty(\Omega)}
   \Norm{\grad u_h^n}_{L^2(\Omega)}
   \Norm{\grad\left( P^n[u^n] - u^n\right)}_{L^2(\Omega)}.
\]
Using again the uniform boundedness of $u_h^n$ in $H^1_0(\Omega)$ provided by \eqref{est:p} of Lemma \ref{lem:Poisson}, combined with the $H^1$-stability \eqref{H1stable1},
we conclude.
\endproof

\subsection{$\Gamma$-convergence of the semi-discrete energy functionals and construction of global energy minimizers}

The main result of this section is the $\Gamma$-convergence of the sequence of the semi-discrete energy functionals \eqref{def:EEh} in the appropriate topology.

\begin{theorem}\label{thm:GammaConv}
Let $(\cc^n)_{n\in\N}$ be a sequence of nonnegative conductivities $\cc^n \in Y_0^n$ such that $sup_{n\in\N} \Norm{c^n}_{L^\infty(\Omega)} <\infty$
and $\cc^n\wto \cc$ weakly-$\ast$ in $L^\infty(\Omega)$. Then
\(  \label{eq:liminf}
   \E[\cc] \leq \liminf_{n\to\infty} \E^n[\cc^n].
\)
Moreover, for any $\cc\in L^\infty(\Omega)$ there exists a sequence $(\cc^n)_{n\in\N} \subset L^\infty(\Omega)$ converging to $\cc$ in the norm topology of $L^q$, $q<\infty$, such that
\(  \label{eq:limsup}
   \E[\cc] \geq \limsup_{n\to\infty} \E^n[\cc^n].
\)
\end{theorem}

\begproof
Statement \eqref{eq:liminf} is a direct consequence of
the weak lower semicontinuity \eqref{eq:wlsc} of the energy functional $\E$, combined with Lemma \ref{lem:Enn}.
Indeed, for any $\cc^n\wto \cc$, we have
\[
   \E[\cc] &\leq& \liminf_{n\to+\infty} \E[\cc^n] \\
   &\leq& \liminf_{n\to\infty} \E^n[\cc^n] + \lim_{n\to\infty} ( \E[\cc^n] - \E^n[\cc^n] ) \\
   &=& \liminf_{n\to\infty} \E^n[\cc^n].
\]

To prove \eqref{eq:limsup}, let us fix $\cc\in L^\infty(\Omega)$ and construct $\cc^n\in Y^n_0$ as the sequence of piecewise constant approximates of $\cc$. By standard approximation theory we then have $\cc^n\to \cc$ in the norm topology of $L^q(\Omega)$ for any $q<\infty$. 
Then the strong continuity of the energy functional $\E$ with respect to the $L^2$-topology established in Lemma \ref{lem:strong} gives
\[
   \E[\cc] = \lim_{n\to+\infty} \E[\cc^n].
\]
Moreover, Lemma \ref{lem:ineqE} gives $\E[\cc^n] \geq \E^n[\cc^n]$, which directly implies \eqref{eq:limsup}.
\endproof

With Theorem \ref{thm:GammaConv} it is straightforward to construct global minimizers of the continuum energy \eqref{energy} as limits of sequences of minimizers of the discrete energy functional \eqref{def:en:disc:resc} on equilateral triangulations.

\begin{theorem}
Let $\CC^n \in \R^{|\Eset^n|}$, $n\in\N$, be a sequence of global minimizers of the rescaled discrete energy functionals $\bEn^n$ given by 
\eqref{def:en:disc:resc}.
Assume that 
$sup_{n\in\N} \Norm{\CC^n}_{\ell^\infty} <\infty$.
Then, up to an eventual selection of a subsequence,
$\cc^n:=\Zh[\CC^n]$ converges weakly-*
to $\cc\in L^\infty(\infty)$.
Under the additional assumption that $\cc\in C(\bar\Omega)$, $\cc$ is a global minimizer of the continuum energy functional \eqref{energy}.
\end{theorem}

\begproof
The uniform boundedness of $\Norm{\CC^n}_{\ell^\infty}$ obviously implies that $\cc^n:=\Zh[\CC^n]$ is uniformly bounded in $L^\infty(\Omega)$ and (an eventual subsequence) converges weakly-* to some $\cc\in L^\infty(\infty)$. Then statement \eqref{eq:liminf} of Theorem \ref{thm:GammaConv} gives
\[
   \E[\cc] \leq \liminf_{n\to+\infty} \E^n[\cc^n].
\]
We claim that if $\cc\in C(\bar\Omega)$, then it is a global minimizer of the energy functional \eqref{energy}.
For contradiction, let us assume that there exists
$\tilde\cc\in C(\bar\Omega)$ such that
\(  \label{eq:contra}
   \E[\tilde\cc] < \E[\cc].
\)
By averaging $\tilde\cc$ over the diamonds $\diamond_{ij}^n$, we construct a sequence $(\widetilde\CC^n)_{n\in\N}$ such that $\tilde\cc^n := \Qh[\widetilde\CC^n]$ converges to $\tilde\cc$ in the norm topology of $L^\gamma(\Omega)$. Moreover, since $\tilde\cc\in C(\bar\Omega)$, the sequence $\widetilde\CC^n$ satisfies the property \eqref{eq:coro}.
Corollary \ref{lem:epsn} and the statement \eqref{eq:limsup} of Theorem \ref{thm:GammaConv} gives then
\[
   \limsup_{n\to\infty}\, \bEn^n[\widetilde\CC^n] = \limsup_{n\to\infty}\, \E^n[\tilde\cc^n]
   \leq \E[\tilde\cc].
\]
On the other hand, 
we have
\[
   \liminf_{n\to\infty}\, \bEn^n[\CC^n] = \liminf_{n\to\infty}\, \E^n[\cc^n]
   \geq \E[\cc].
\]
Since by construction we have $\bEn[\CC^n] \leq \bEn[\widetilde\CC^n]$ for all $n\in\N$, we finally obtain
\[
   \E[\cc] \leq \liminf_{n\to\infty}\, \bEn^n[\CC^n]
   \leq
   \limsup_{n\to\infty}\, \bEn^n[\widetilde\CC^n]
   \leq \E[\tilde\cc],
\]
which is a contradiction to \eqref{eq:contra}.
\endproof

Finally, let us note that by the standard theory of convex gradient flows on Hilbert spaces, see, e.g., \cite{AGS, Santambrogio}, we can construct unique weak solutions of 
the system \eqref{eq:Poisson}--\eqref{eq:A}.
We refer to \cite[Section 3.1]{HMP22} for further details.

\begin{proposition}
Let $r>0$, $S\in L^2(\Omega)$ satisfying \eqref{ass:S}
and $\cc^0 \in L^2(\Omega)$ nonnegative almost everywhere in $\Omega$.
Then the problem \eqref{eq:Poisson}--\eqref{IC_0}
admits a unique global weak solution $\cc \in H^1((0,\infty); L^2(\Omega))$,
with $\En_0[\cc]\in L^\infty(0,\infty)$ and satisfying the energy inequality
\[
   \E[\cc(t)] + \int_0^t \int_\Omega \left| \part{\cc}{t}(s,x) \right|^2 \d x\d s \leq \E[\cc^0].
\]
\end{proposition}

\section{Numerical results}\label{sec:numerics}

In this Section, we briefly discuss some of the details of our numerical implementation. We first introduce the finite element spaces,  the semidiscrete nonlinear equations and their time advancing scheme, and finally discuss the linearized equations. We conclude with numerical experiments in Section \ref{sec:sequence_formation} and \ref{subsec:pLaplacian}. For additional details, see \cite{HMPZ}. 

Here we use finite element spaces consisting of piecewise constant elements for the conductivity and bilinear elements for the potential in two dimensions, i.e.,
\begin{align*}
\Cspace_h &= \big\{ c_h \in L^2(\Omega) ~|~ {c_h}_{|K} \text{ constant, } \forall K \in \Omega\big\},\\
\Pspace_h &= \big\{ u_h \in H_{\Gamma_D}^1(\Omega) ~|~ u_h \in C^0(\Omega),\, u_h=g_D \text{ on } \partial\Gamma_D, \text{ and } {u_h}_{|K} \text{ bilinear, } \forall K \in \Omega\big\},
\end{align*}
where $K$ is a quadrilateral cell of the tessellated domain $\Omega$ and $\Gamma_D\cup \Gamma_N = \partial \Omega$, with $\Gamma_D$ (resp. $\Gamma_N$) the subset of the boundary where we may consider non-homogeneous essential (resp. natural) boundary conditions with datum $g_D$ (resp $g_N$). Clearly, for the case of non-homogeneous natural boundary conditions, the model equations are complemented with
\[
c \frac{\partial u}{\partial n} = g_N\quad \mbox{ on } \Gamma_N.
\]

For the sake of numerical stability, we may use a regularized metabolic cost
\[
(c^2 + \eps)^{\gamma/2},\quad \eps \ge 0,
\]
that reduces to $c^{\gamma}$ when $\eps=0$.
After a simple rescaling to ensure the symmetry of the linearized equations, the variational formulation of the discrete problem reads: given $S \in L^2(\Omega)$, find $u_h \in \Pspace_h$, $c_h \in \Cspace_h$ such that
\begin{align*}
&\frac{1}{2}\int_\Omega \left( \dd{c_h}{t} \, b_h - |\nabla u_h|^2 \, b_h + {\nu} \,(c_h^2 + \eps)^{(\gamma - 2)/2} c_h \,b_h \right) \dx = 0, &\forall\, b_h \in \Cspace_h,\\
&-\int_\Omega(c_h + r)\nabla u_h \cdot \nabla q_h ~ \dx = - \int_\Omega S q_h ~ \dx  -\int_{\Gamma_N} g_N \,q_h ~ \ds, &\forall\, q_h \in \Pspace_h,
\end{align*}
which leads to the set of differential-algebraic equations
\begin{equation*} 
\begin{aligned} %
&\frac{1}{2}\int_\Omega \left( \dd{c_h}{t} \, \Ctest_i - |\nabla u_h|^2 \, \Ctest_i + {\nu} \,(c_h^2 + \eps)^{(\gamma - 2)/2} c_h \, \Ctest_i \right) \dx = 0,\\
&-\int_\Omega(c_h + r)\nabla u_h \cdot \nabla \ptest_i ~ \dx = - \int_\Omega S \ptest_i ~ \dx -\int_{\Gamma_N} g_N \,\ptest_i ~ \ds,
\end{aligned}
\end{equation*}
where $\ptest_i$ (resp. $\Ctest_i$) are test functions for $\Pspace_h$ (resp. $\Cspace_h$).

To advance in time, we use the backward Euler scheme, which allows us to preserve the positivity of the conductivity. Dropping the $h$ subscript, we first compute the initial potential $u_0$ by solving
\begin{equation*} 
\int_\Omega(c_0 + r)\nabla u_0 \cdot \nabla \ptest_i ~ \dx =  \int_\Omega S\, \ptest_i ~ \dx + \int_{\Gamma_N} g_N \,\ptest_i \ds,
\end{equation*}
where $c_0$ is the given initial condition for the conductivity.
Denoting by $c_n$ and $u_n$ the finite element approximations for the conductivity and potential at step $n$, we then solve the equations
\begin{equation*} 
\begin{aligned} %
&\frac{1}{2}\int_\Omega \left( \frac{c_{n+1} - c_{n}}{\delta t_n} \, \Ctest_i - |\nabla u_{n+1}|^2 \, \Ctest_i + {\nu} \,(c_{n+1}^2 + \eps)^{(\gamma - 2)/2} c_{n+1} \, \Ctest_i \right) \dx = 0,\\
&-\int_\Omega(c_{n+1} + r)\nabla u_{n+1} \cdot \nabla \ptest_i ~ \dx = - \int_\Omega S \ptest_i ~ \dx - \int_{\Gamma_N} g_N \,\ptest_i \ds,\\\end{aligned}
\end{equation*}
where $\delta t_n$ is the time step.

This nonlinear system of equations is solved using the inexact Newton method \cite{NocedalWright}. The symmetric indefinite Jacobian matrix at a given linearization point is of the form
\[
J = \begin{bmatrix}
A & B^T \\
B & -C
\end{bmatrix},
\]
where
\begin{align*}
&A_{ij} =  \frac{1}{2}\int_\Omega \bigg(\frac{1}{\delta t_n}  + \alpha(c_{n,k})+ \beta(c_{n,k})  \bigg) \ctrial_i \ctrial_j \dx,\\
&B_{ij} = -\int_\Omega \ctrial_j \left(\nabla u_{n,k} \cdot \nabla \ptrial_i\right)\dx,\\
&C_{ij} = \int_\Omega  (c_h + r ) \nabla \ptrial_i \cdot \nabla \ptrial_j \dx,
\end{align*}
with $c_{n,k}$ (resp. $u_{n,k}$) is the linearization point at the $k$th Newton step for the conductivity (resp. potential), and the function coefficients $\alpha$ and $\beta$ are given by
\begin{equation*} 
\alpha(c) = \nu(c^2 + \eps)^{(\gamma - 2)/2},\quad \beta(c) = \nu\,(\gamma - 2)(c^2 + \eps)^\frac{\gamma - 4}{2}c^2.
\end{equation*}
Note that, when $\gamma \ge 1$ and $\eps \ge 0$, $A$ is a symmetric positive definite matrix because
\[
\frac{1}{\delta t_n} + \alpha(c) + \beta(c) \ge \nu(c^2 + \eps)^{(\gamma - 4)/2}(\eps + (\gamma - 1)c^2) \ge 0.
\]
When $\gamma < 1$, the positive definiteness of $A$ can not be guaranteed unless we make assumptions on the time step.

The matrix $J$ admits the exact factorization
\begin{equation}\label{eq:schurPre}
\begin{bmatrix}
A & B^T \\
B & -C
\end{bmatrix} =
\begin{bmatrix} 
A & 0\\
B & I
\end{bmatrix}\,
\begin{bmatrix} 
A^{-1}& 0\\
0 & -S
\end{bmatrix} \,
\begin{bmatrix} 
A & B^T\\
0 & I 
\end{bmatrix},
\quad S = C + B A^{-1} B^T,
\end{equation}
where the Schur complement $S$ can be explicitly assembled since $A$ is a block-diagonal matrix with each block associated with a single element. The inverse factorization requires solving the inexpensive diagonal problem $A$ twice, and the Schur complement system $S$ only once. Numerically, we observed that $S$ is a symmetric positive semi-definite matrix representing a perturbed scalar diffusion problem having the same kernel as the discretized Poisson problem $C$. Instead of solving the Schur complement system exactly, we can thus use an algebraic multigrid (AMG) preconditioner, which guarantees robustness and has an overall cost that is linear in the matrix size. 

The code used to conduct the numerical experiments described in this Section is publicly available as a tutorial of the Portable and Extensible Toolkit for Scientific Computing (PETSc, \cite{petsc})\footnote{\url{https://gitlab.com/petsc/petsc/-/blob/main/src/ts/tutorials/ex30.c}}. The management of unstructured grids and the implementation of finite elements are based on the {\tt DMPLEX} infrastructure \cite{dmplex}. Time integration is performed using the {\tt TS} module \cite{abhyankar2018petsc}; adaptive time stepping is performed using digital filtering-based methods \cite{soderlind2003} combined with the approximation of the local truncation error using extrapolation. The nonlinear systems of equations arising from the time-advancing schemes are solved with a tight absolute tolerance of 1.E-10; Jacobian linear systems are solved using the right-preconditioned GMRES method \cite{saad1986}, and the Schur complement system is approximatively solved using the application of an AMG preconditioner based on smoothed aggregation \cite{adams2003}. Linear system tolerances are dynamically adjusted using the well-known Eisenstat and Walker trick, which prevents over-solving in the early steps of the Newton process and tightens the accuracy as convergence is approached \cite{ew96}. The nonlinear problems can sometimes be quite difficult to solve; as a safeguard, we cap the maximum number of nonlinear iterations to 10 and shorten the time step if the Newton process stagnates and does not reach convergence within the maximum number of iterations.

\subsection{Network formation}
\label{sec:sequence_formation}

\begin{figure}[htbp]
\centering
\begin{tabular}{c c c}
\includegraphics[width=0.3\textwidth,trim={50 10 80 0},clip]{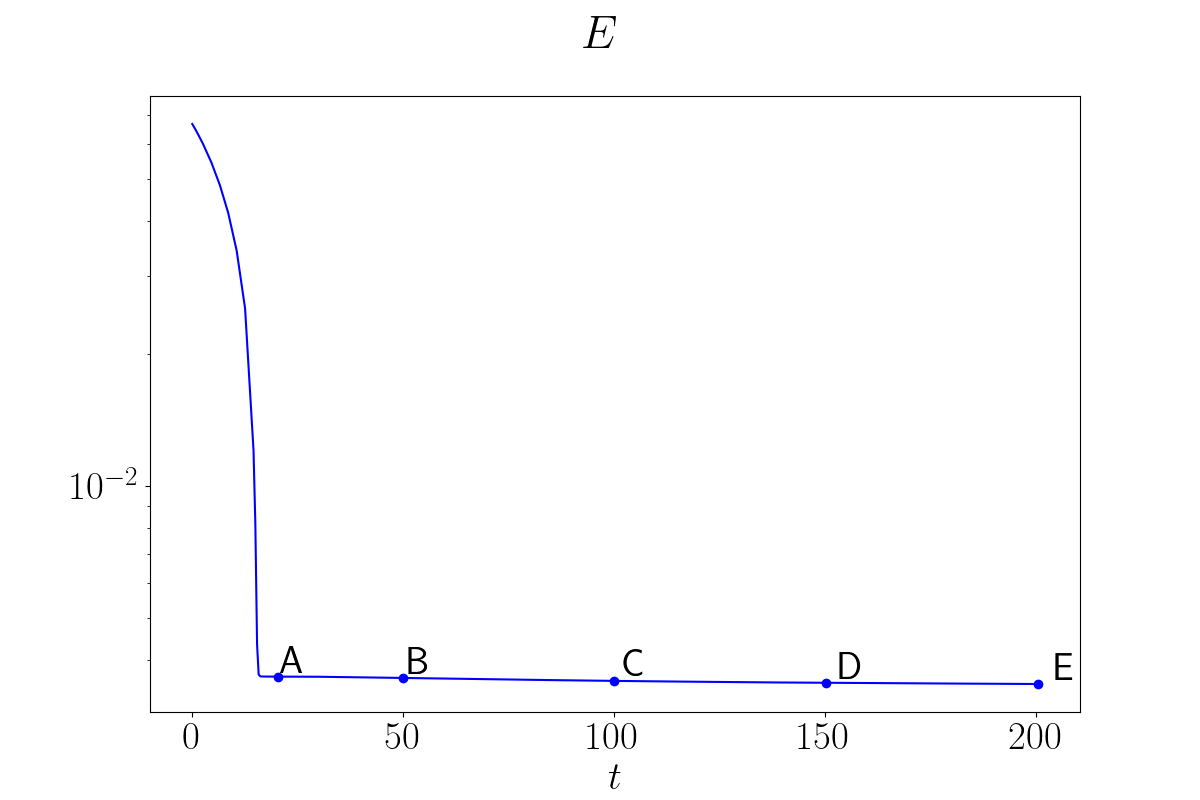} & \includegraphics[width=0.3\textwidth,trim={50 10 80 0},clip]{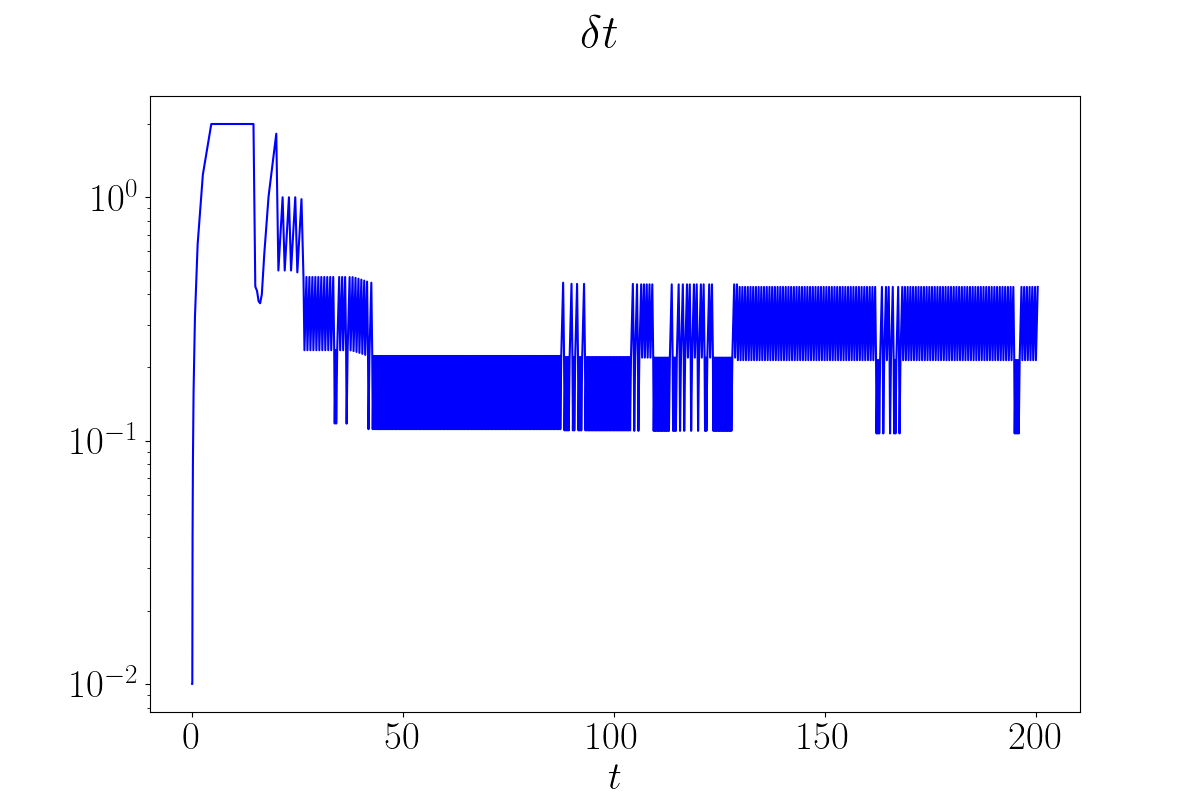}
\includegraphics[width=0.3\textwidth,trim={50 10 80 0},clip]{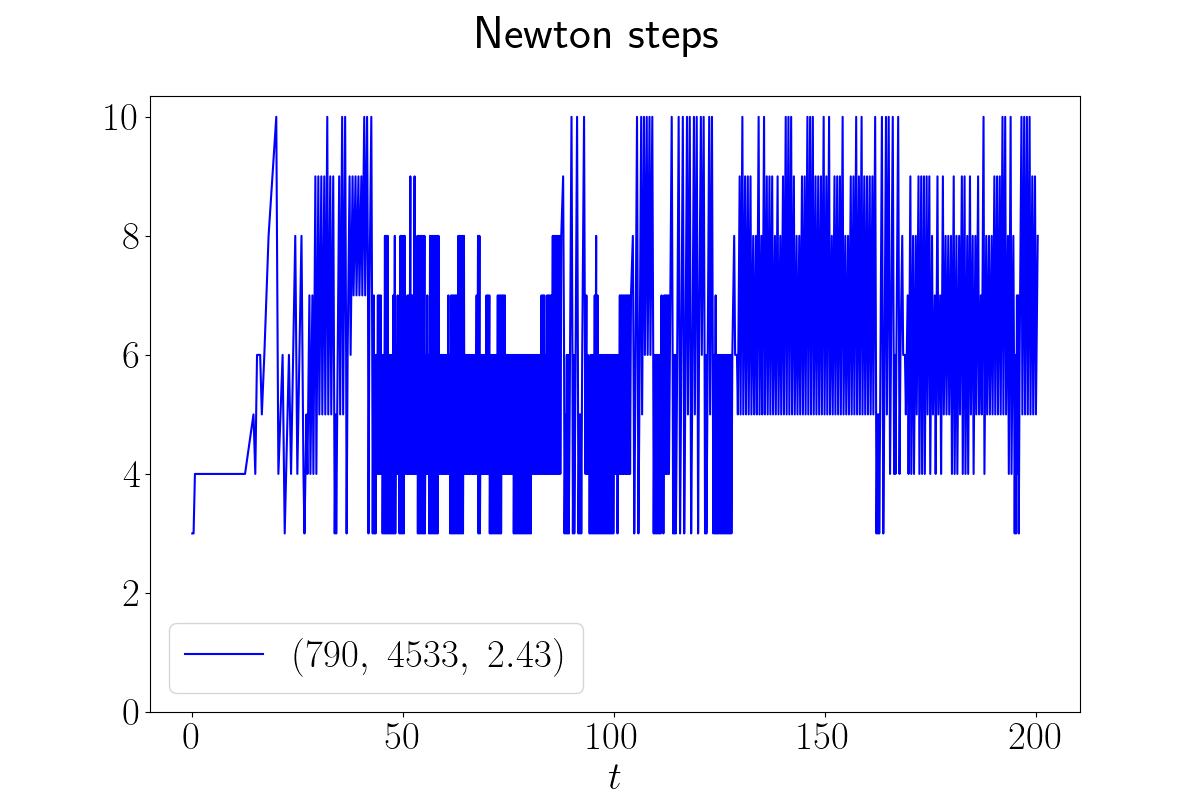}
\end{tabular}
\caption{Network formation: Energy (left panel), time step (central), and number of nonlinear iterations as a function of simulation time for 512x512 mesh. The legend in the right panel reports in parentheses the total number of time steps, the total number of nonlinear steps, and the average number of Krylov iterations per Newton step.
}
\label{fig:box_sequence_logs}
\end{figure}

\begin{figure}[htbp]
\centering
\begin{tabular}{c c c c c}
A & B & C & D & E\\
\includegraphics[width=0.17\textwidth,trim={110 70 110 70},clip]{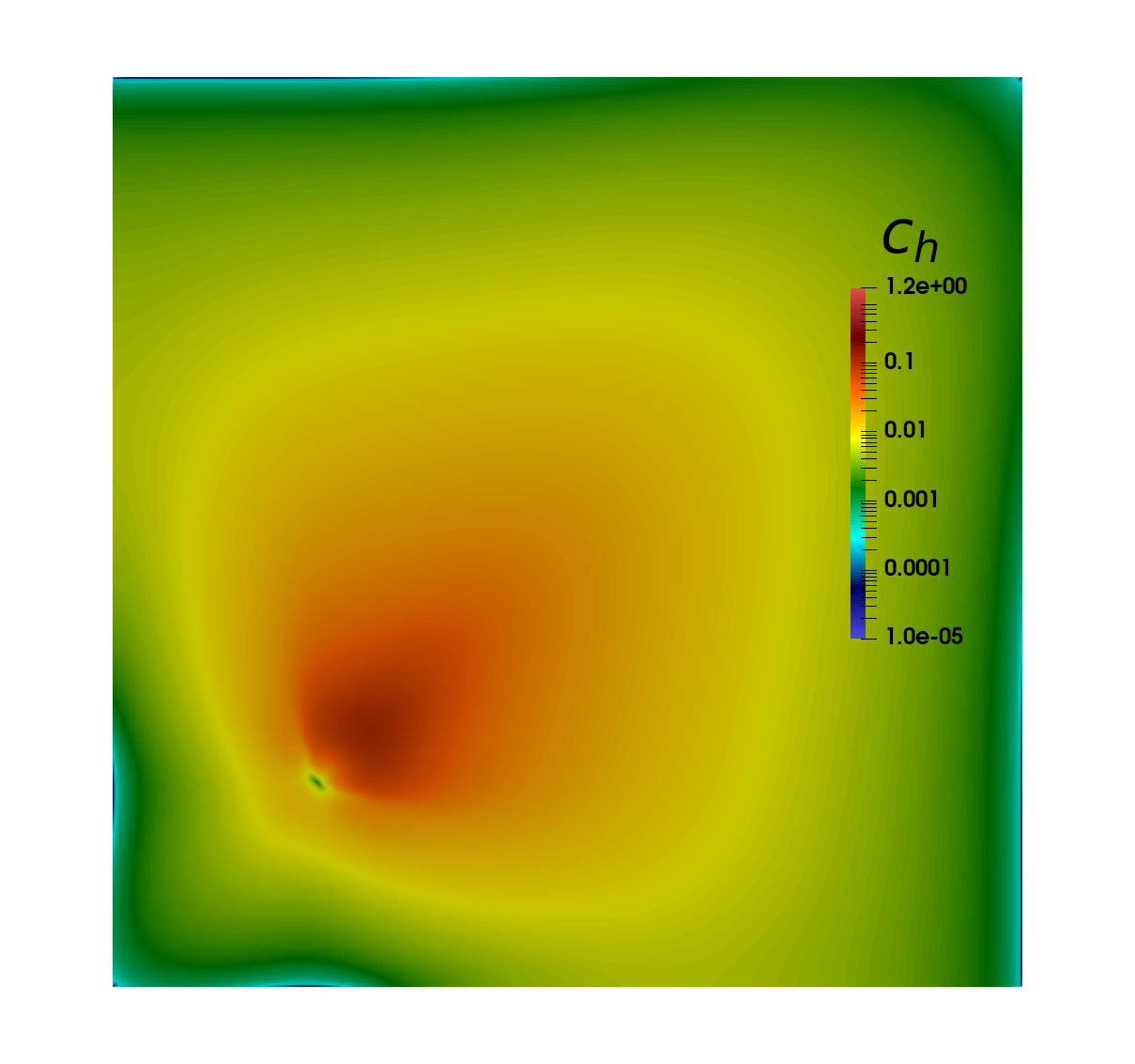} & 
\includegraphics[width=0.17\textwidth,trim={110 70 110 70},clip]{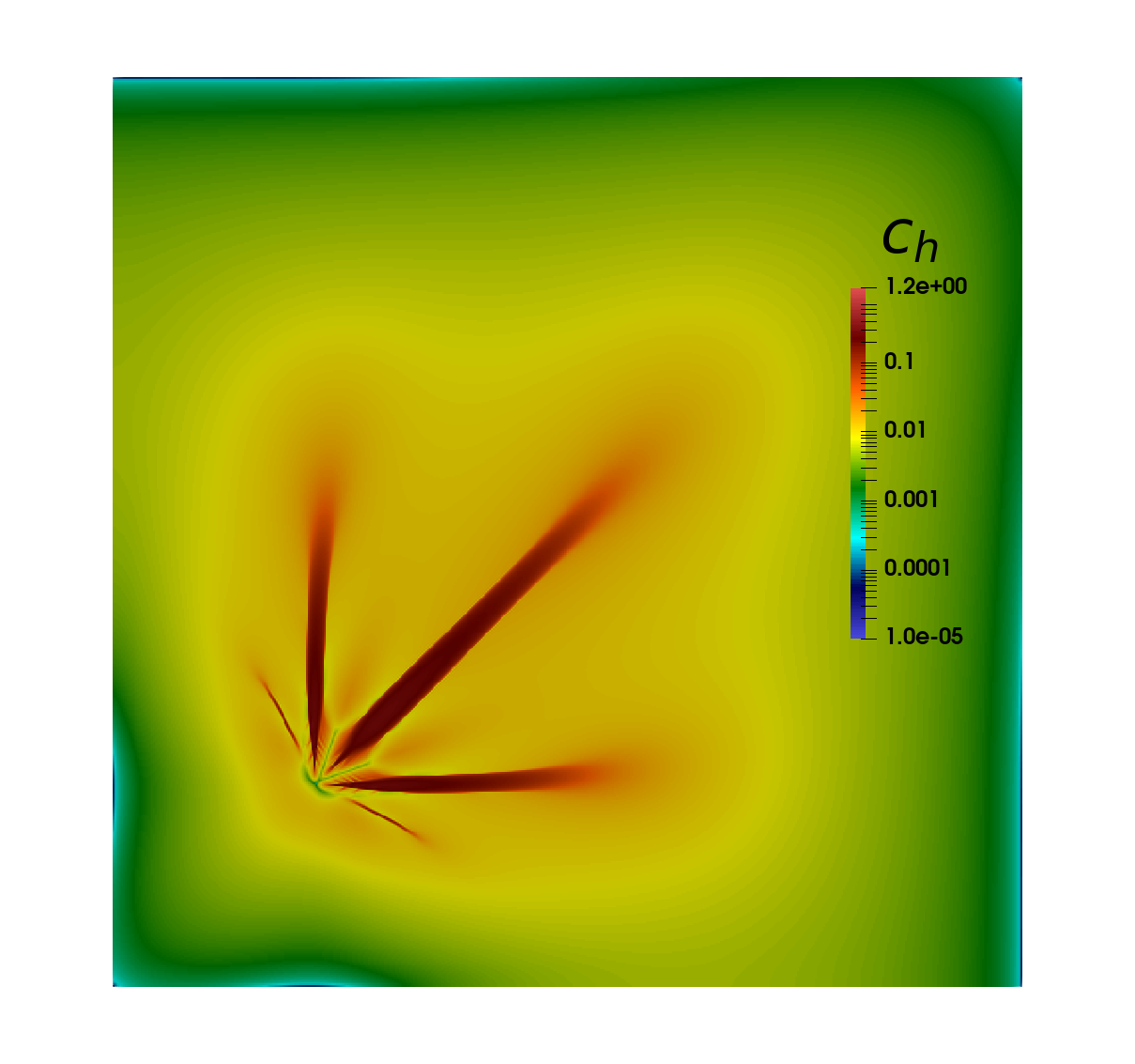} & 
\includegraphics[width=0.17\textwidth,trim={110 70 110 70},clip]{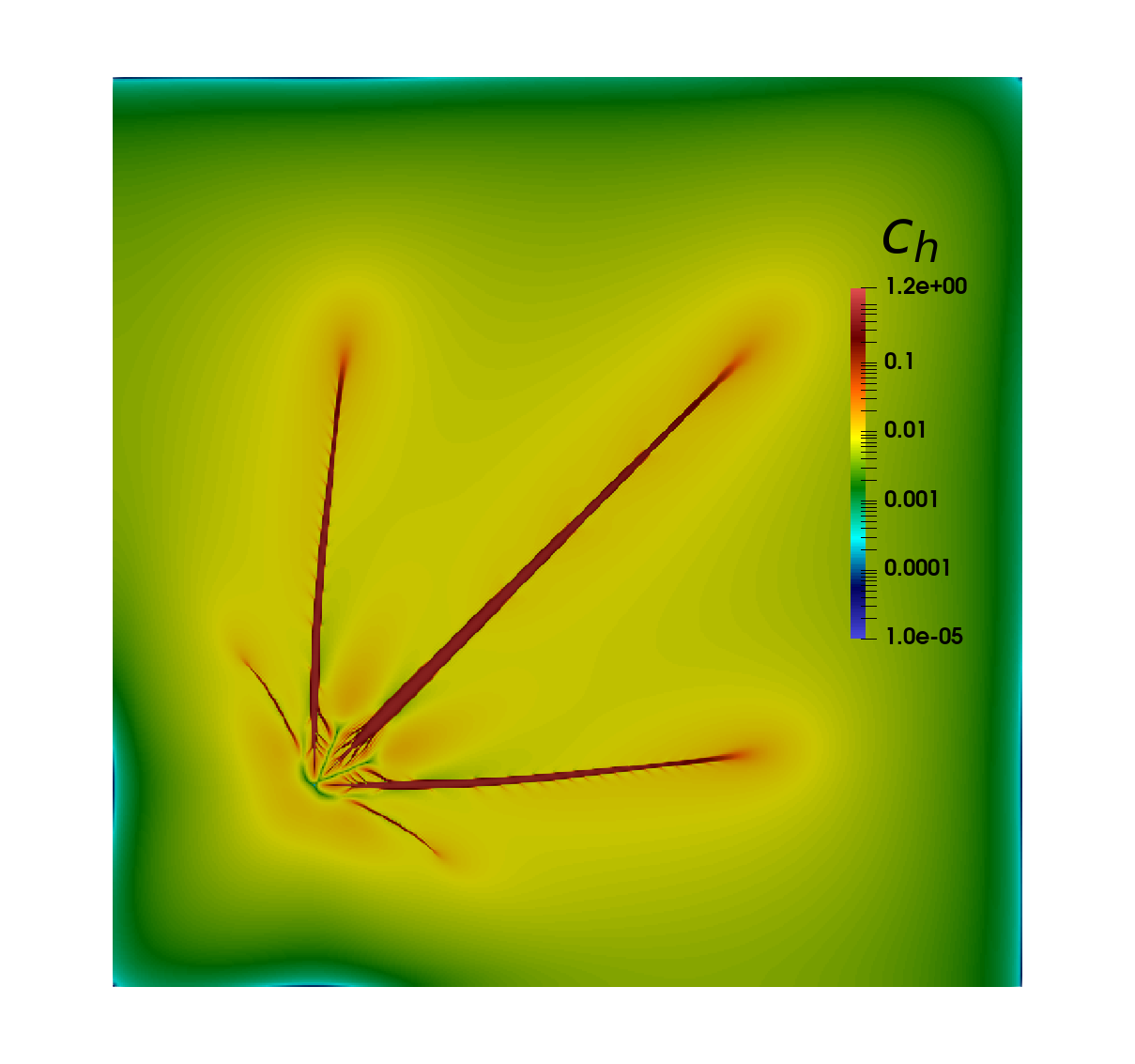} & 
\includegraphics[width=0.17\textwidth,trim={110 70 110 70},clip]{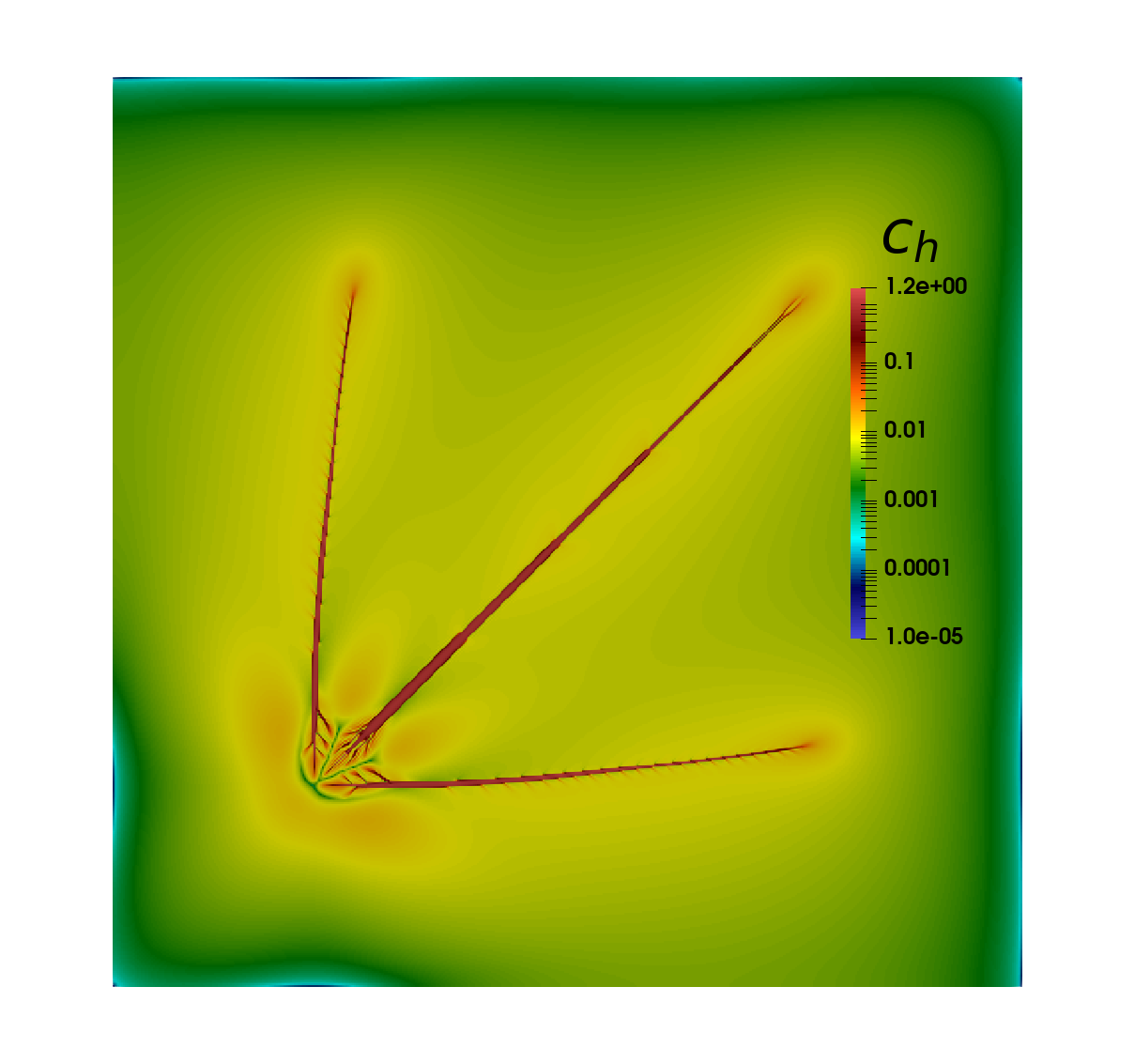} & 
\includegraphics[width=0.17\textwidth,trim={110 70 110 70},clip]{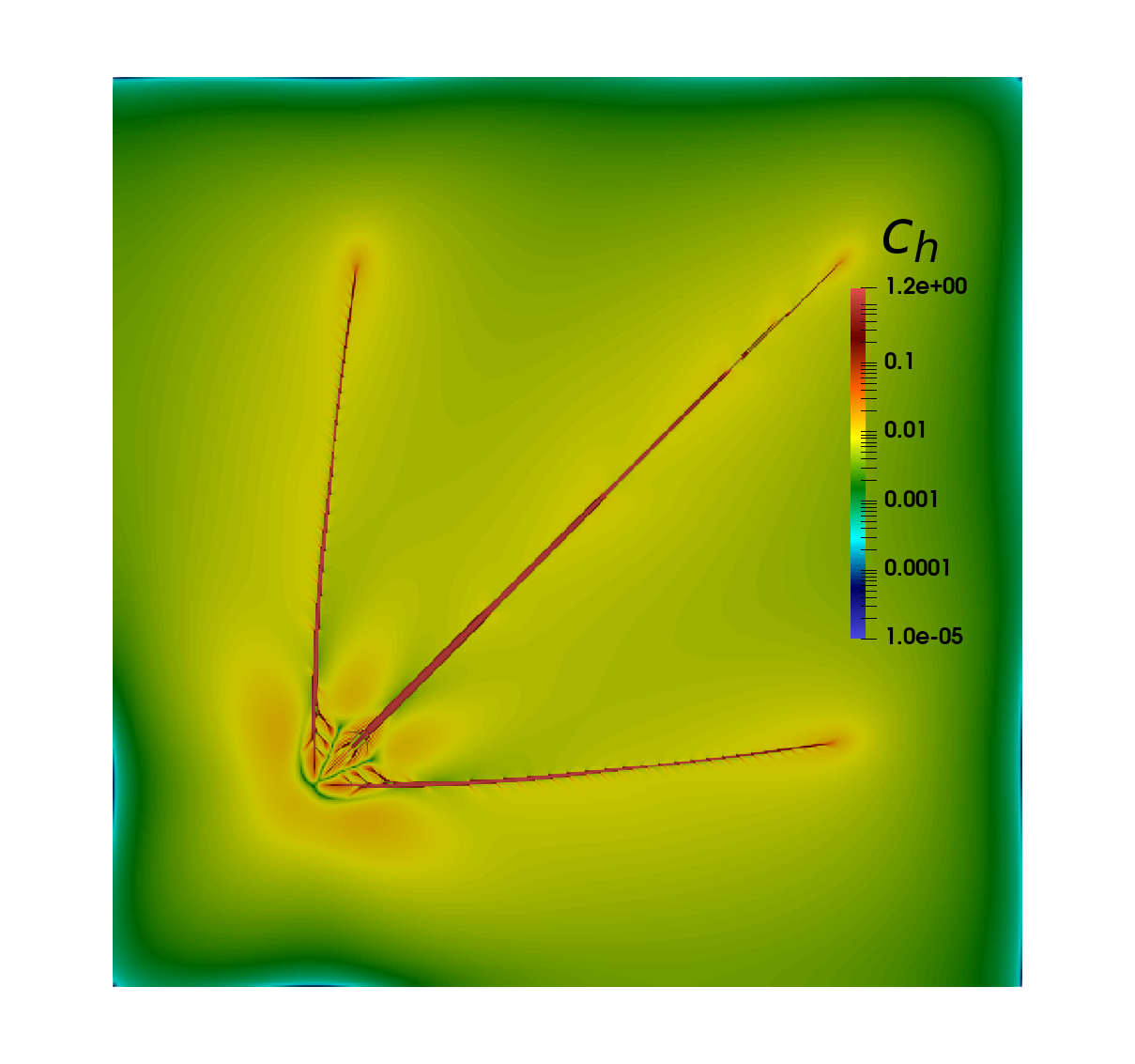} 
\end{tabular}
\caption{Network formation: Conductivity in logarithmic scale at selected time instances (see Figure \ref{fig:box_sequence_logs}).
}
\label{fig:box_sequence_imgs}
\end{figure}

In the first numerical experiment, we prove that the scalar model is able to simulate a network formation process. For this experiment, we consider a unit square domain $\Omega=[0,1]^2$, homogeneous natural boundary conditions, and the model parameters $r=0.0001$, $\eps=0.001$, $\gamma=0.75$, and $\nu=0.05$.  The domain is discretized with a structured quadrilateral mesh with a 512x512 cell arrangement; the total number of degrees of freedom is 263 thousand. The initial conductivity field is the constant function  $c_0=1$, the simulation is run up to time $t=200$, and we use the Gaussian source
\begin{equation*}
S(x) = S_0(x) - \int_\Omega S_0(x) ~\dx, \quad S_0(x) = e^{-500\Norm{x-x_0}^2},
\end{equation*}
with $x_0=(0.25,0.25)$. Figure \ref{fig:box_sequence_logs} presents the system energy $E$ (left panel), the time step size $\delta t_n$ (center), and the number of nonlinear iterations (right) as functions of the simulated time. The legend in the right panel reports the total number of time steps, the total number of Newton steps, and the average number of Krylov iterations per Newton step to provide an overview of the algorithmic costs of the solver. Snapshots of the conductivity field $c_h$ (in logarithmic scale) taken at selected times indicated by the text boxes in Figure \ref{fig:box_sequence_logs}, are shown in Figure \ref{fig:box_sequence_imgs}.

The system energy decreases monotonically throughout the simulation. In the initial phase, the energy decays rapidly as the conductivity field diffuses across the domain (see panel A of Figure \ref{fig:box_sequence_imgs}), allowing relatively large time steps. At this stage, the dynamics are dominated by the minimization of the term $(c + r)|\nabla u|^2$. In the subsequent phase, network formation begins (see panels B and C), driven by the interplay between this positive term and the negative contribution from the metabolic term. During this process, energy variations become smaller, and the solver selects smaller time steps. In the final phase, the energy decreases more slowly toward its steady state, with minimal changes in the conductivity field (see panels D and E).
Overall, the problem exhibits a mild nonlinearity. The backward Euler time advancing scheme encounters a few nonlinear solver failures during the later stages of network formation. Nevertheless, the average number of Krylov iterations per Newton step remains very low, confirming the effectiveness and robustness of the linear preconditioning strategy.

\subsection{Solutions of the $p$-Laplacian equation}
\label{subsec:pLaplacian}

We dedicate a substantial portion of the numerical section to the integration of the differential-algebraic equations \eqref{eq:Poisson} and \eqref{eq:A} up to their steady state, in order to compute solutions of the $p$-Laplacian equations
\[
-\nabla  \cdot (|\nabla u|^{p-2} \nabla u) = S, \quad p = \frac{2\gamma}{\gamma - 1} > 2,
\]
with $\gamma > 1$, $r=0$, $\eps=0$, and $\nu=1$. To evaluate the approximation properties of the method, we first apply the method of manufactured solutions and compare the numerical results against analytical solutions. Finally, we present experiments demonstrating the robustness of the method for large values of $p$ and non-convex domains

A priori estimates for the lowest order $H^1$-conforming finite element discretizations of the $p$-Laplacian equations with $p>2$ have been established in \cite{EbmeyerLiu2005} (see also \cite{Diening2007}) in terms of a quasi-norm for convex domains and a sufficiently smooth source term $S$ 
\begin{equation}\label{eq:plap_conv}
|u-u_h|^2_{(u,p)} \leq C h^2 \int_\Omega |\nabla u|^{p-2} |D^2 u|^{2} \dx,
\end{equation}
where the quasi-norm $|\cdot|_{(w,p)}$ is defined as 
\begin{equation}\label{eq:quasinorm}
|v|^2_{(w,p)} = \int_\Omega \big(|\nabla w| + |\nabla v|\big)^{p-2} |\nabla v|^{2} \dx.
\end{equation}
In addition to the quasi-norm, we report convergence rates in terms of the standard Lebesgue $L^{p}(\Omega)$ norm and the Sobolev semi-norm defined for $v \in W^{1,p}(\Omega)$
\[
\|v\|_p = \bigg(\int_\Omega |v|^p \dx \bigg)^{1/p},  \qquad |v|_p = \bigg(\int_\Omega |\nabla v|^p \dx \bigg)^{1/p}.
\]

We consider three test cases from \cite{BarrettLiu1993} with radially symmetric smooth functions on the unit square. The first two tests use
\begin{equation}\label{eq:bl1}
\begin{split}
u(r) &= \frac{p - 1}{\sigma+p}(\sigma + 2)^{\frac{1}{1-p}}(1 - r^{\frac{\sigma+p}{p-1}})\\
S(r) &= r^\sigma,\\
\end{split}
\end{equation}
with $\sigma=0, ~p=4$ (denoted by TC1) and $\sigma=7,~p=4$ (TC2) as given in examples 5.2 and 5.3 in \cite{BarrettLiu1993}. The third test case (denoted by TC3) is also from \cite{BarrettLiu1993} (Example 5.4)
\begin{equation}\label{eq:bl2}
\begin{split}
u(r) &= 
\begin{cases}
0& r < a\\
(r-a)^4& r\geq a\\
\end{cases}\\ 
S(r) &=  
\begin{cases}
0& r < a\\
4^{p-1}(r-a)^{3p-4}(2 -3p + \frac{a}{r})& r\geq a\\
\end{cases}\\
\end{split}
\end{equation}
with $a=0.3$ and $p=4$. In order to showcase the convergence of the method for larger $p$-Laplacian exponents, we complement these tests with a fourth one, denoted by TC4, considering \eqref{eq:bl1} with $\sigma=7$  and $p=20$. For all these tests, we use non-homogeneous essential boundary conditions defined by the analytical solution at the boundary.

\begin{figure}[htbp]
\centering
\begin{tabular}{c c}
TC1 $(p=4, \gamma=2)$ & TC2 $(p=4, \gamma=2)$\\
\includegraphics[width=0.4\textwidth,trim={50 0 70 50 },clip]{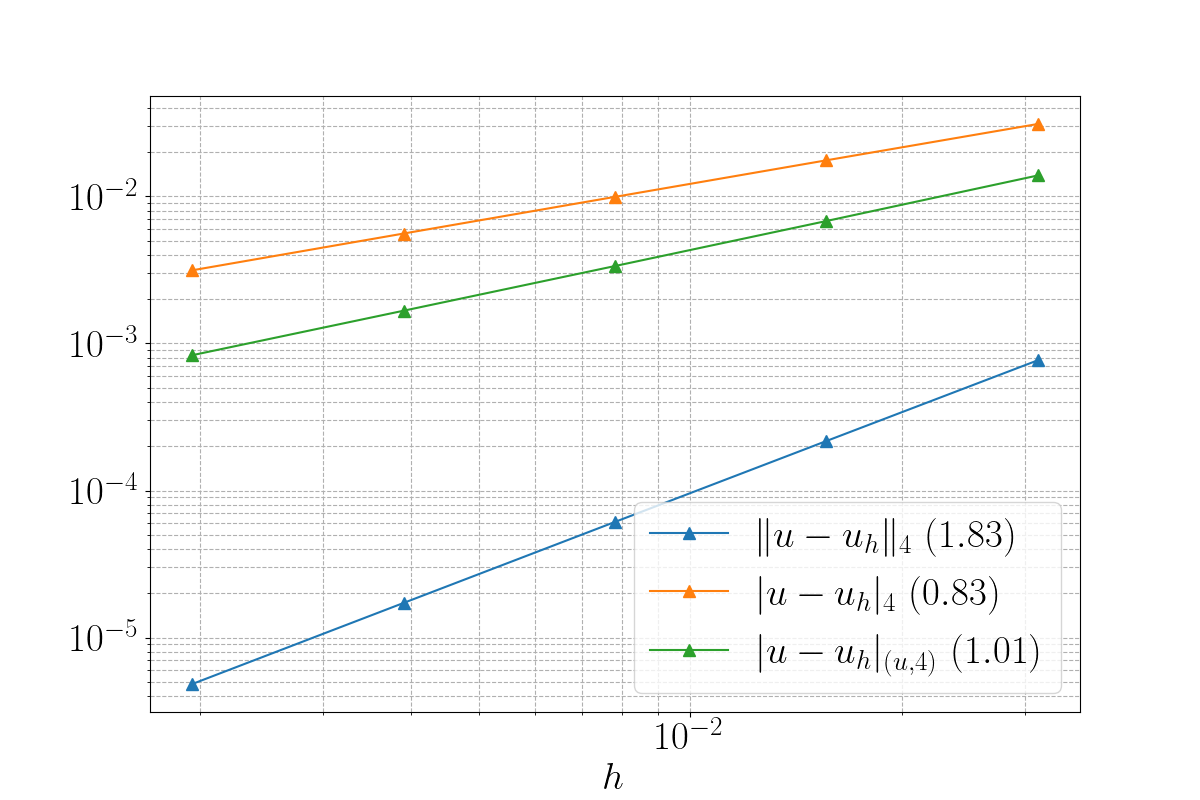} &
\includegraphics[width=0.4\textwidth,trim={50 0 70 50 },clip]{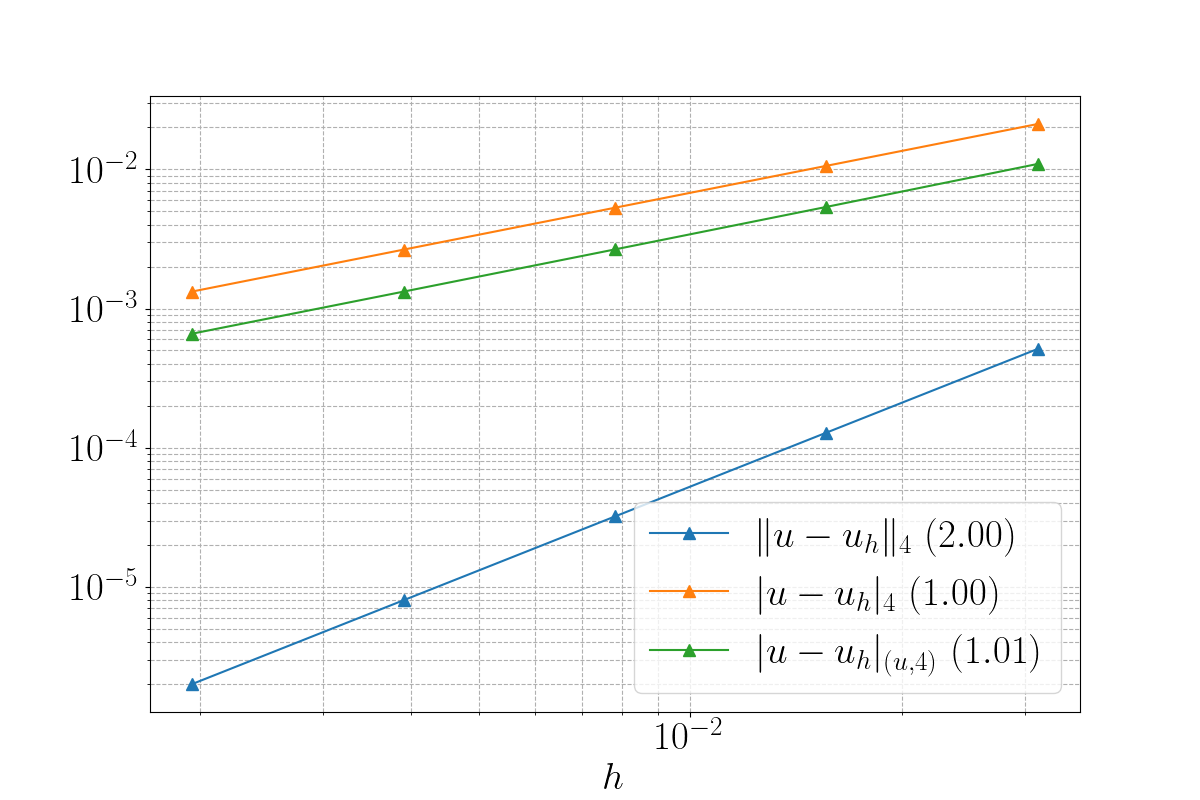} \\
TC3 $(p=4, \gamma=2)$ & TC4 $(p=20, \gamma=10/9)$ \\
\includegraphics[width=0.4\textwidth,trim={50 0 70 50 },clip]{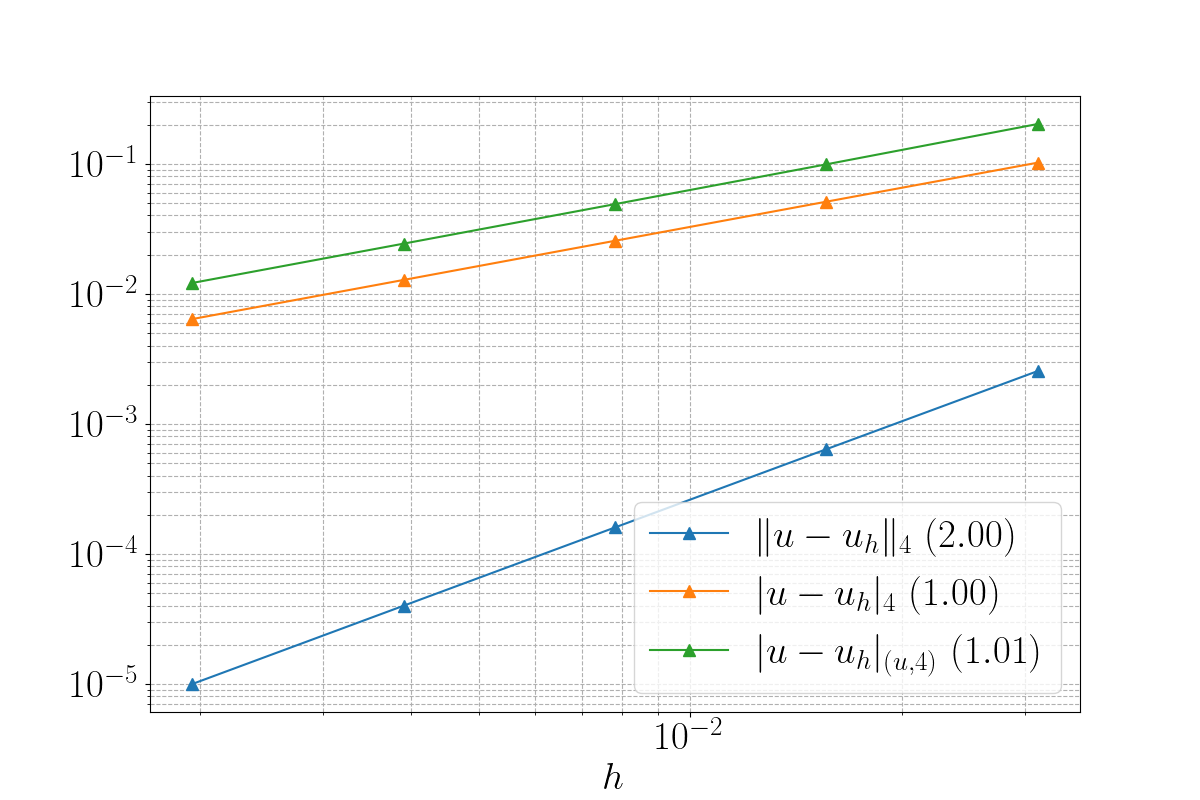} &
\includegraphics[width=0.4\textwidth,trim={50 0 70 50 },clip]{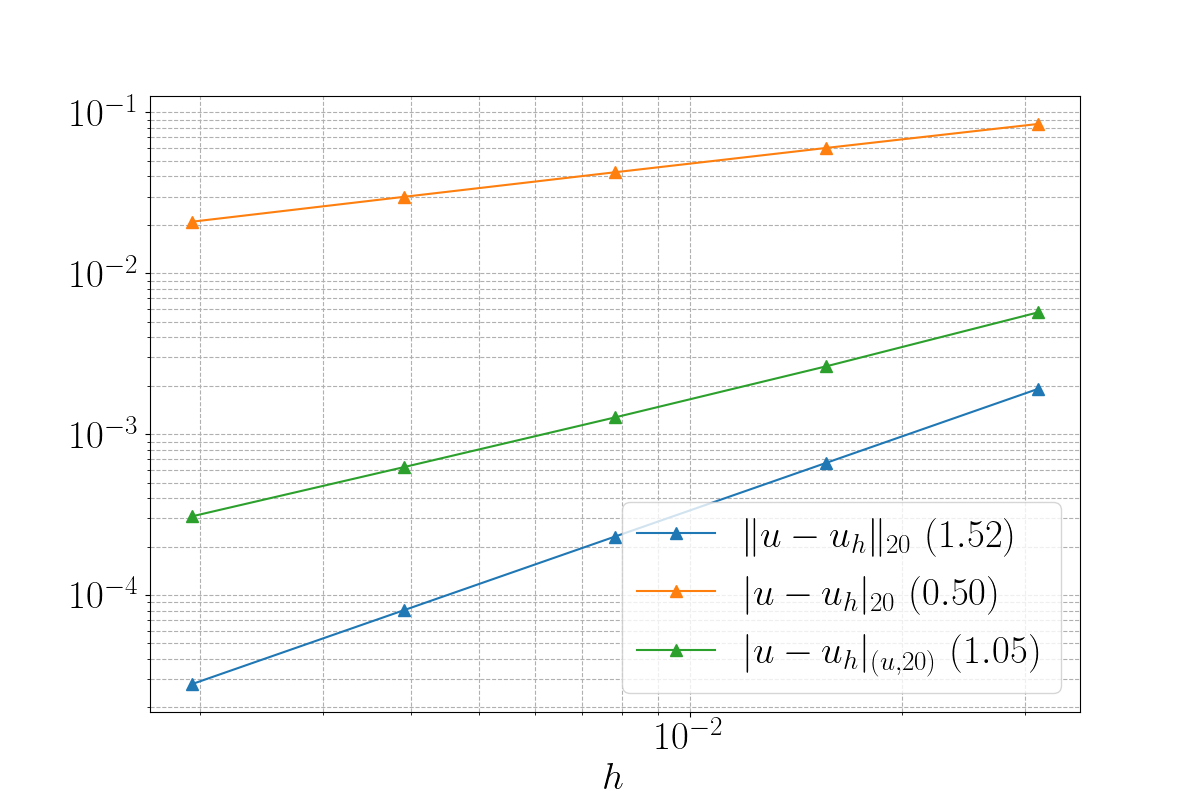} 
\end{tabular}
\caption{MMS errors and convergence rates (in parentheses) for the test cases TC1, TC2, TC3, and TC4 and different error metrics.}
\label{fig:plap_bl_conv}
\end{figure}

We report the convergence rates (in parentheses) in Figure \ref{fig:plap_bl_conv}, considering five levels of uniform refinement starting from a uniform grid of $16 \times 16$ elements. Optimal first-order convergence in the quasi-norm $|\cdot|_{(u,p)}$ is obtained in all the cases. Convergence rates in $L^p(\Omega)$ and $W^{1,p}(\Omega)$ instead depend on the test case. For TC1, we recover the same convergence rates as in \cite{BarrettLiu1993} (see Table 5.2 in the reference), while our convergence rates for TC2 and TC3 are worse than those presented in Tables 5.3 and 5.4 in \cite{BarrettLiu1993}; in particular we always obtain first order convergence in $W^{1,p}(\Omega)$, while the convergence rates reported by Barrett and Liu are approximatively 1.5.

\begin{figure}[htbp]
\centering
\begin{tabular}{c c c c}
& \multicolumn{2}{c}{TC2 $(p=4, \,\gamma=2)$} & \\
\includegraphics[width=0.22\textwidth,trim={50 0 70 0},clip]{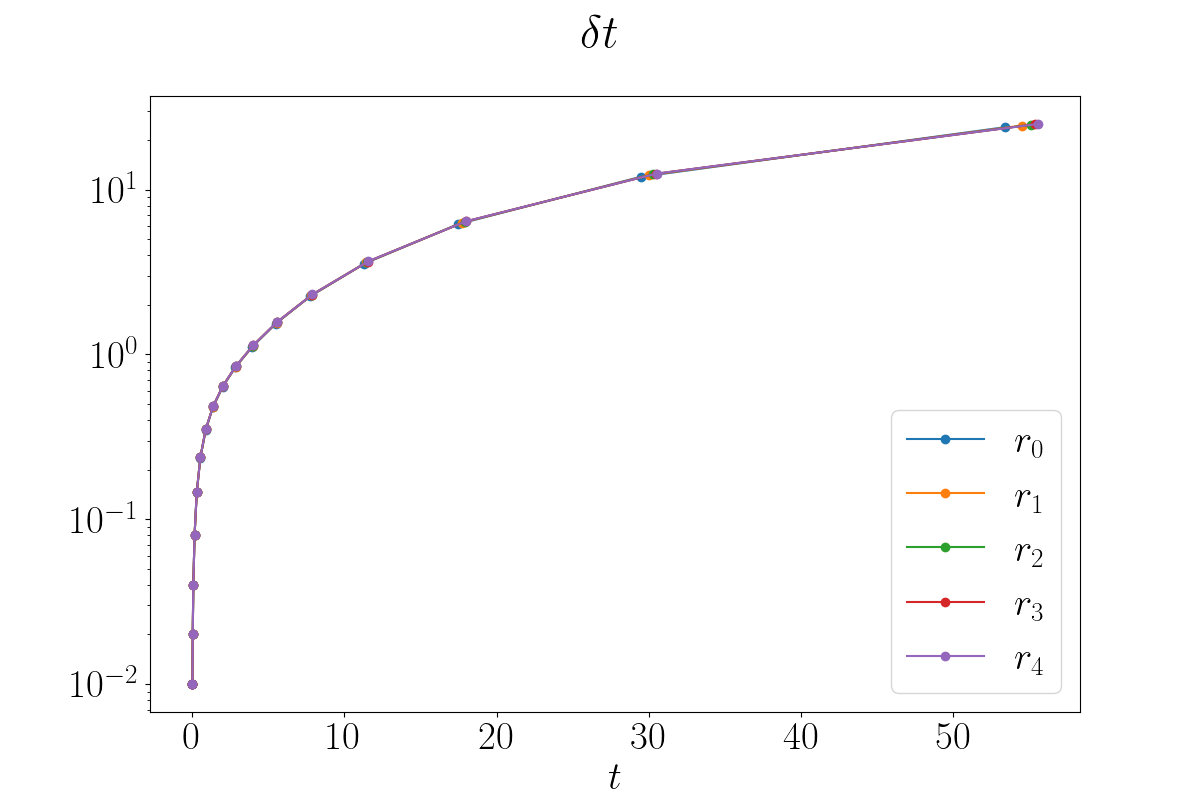} &
\includegraphics[width=0.22\textwidth,trim={50 0 70 0},clip]{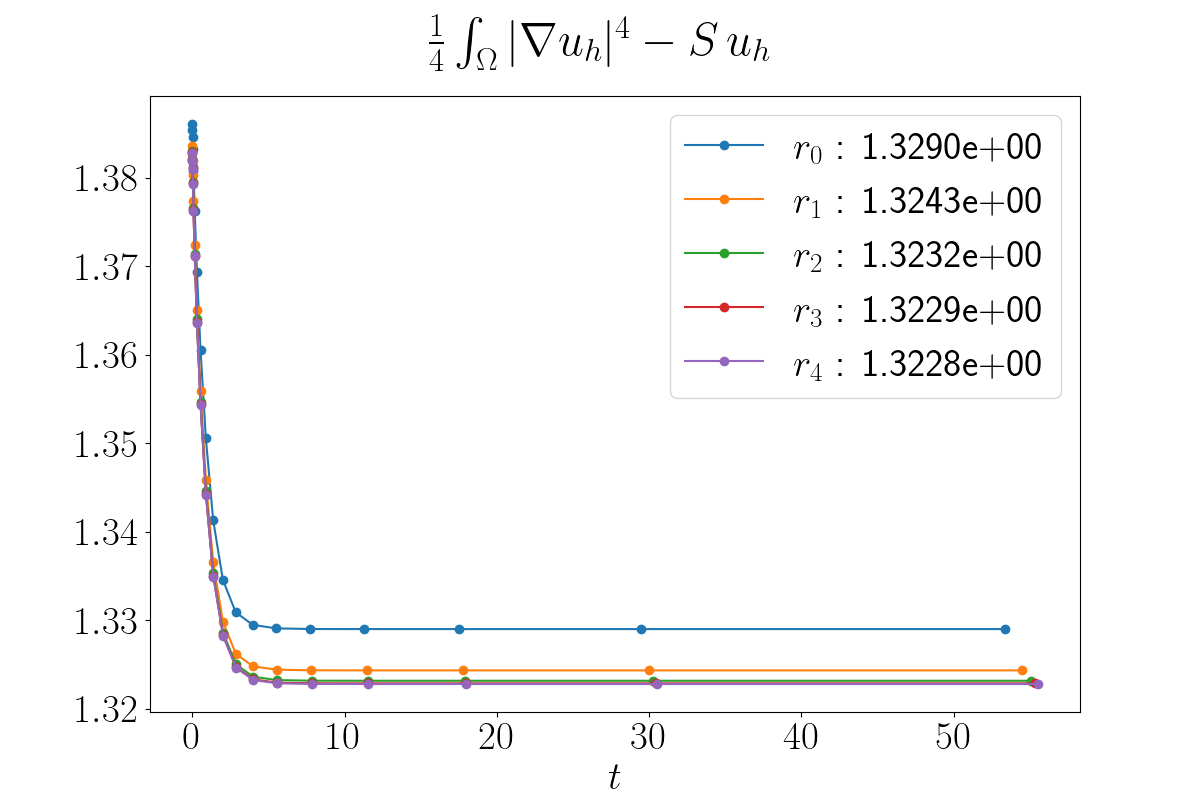} &
\includegraphics[width=0.22\textwidth,trim={50 0 70 0},clip]{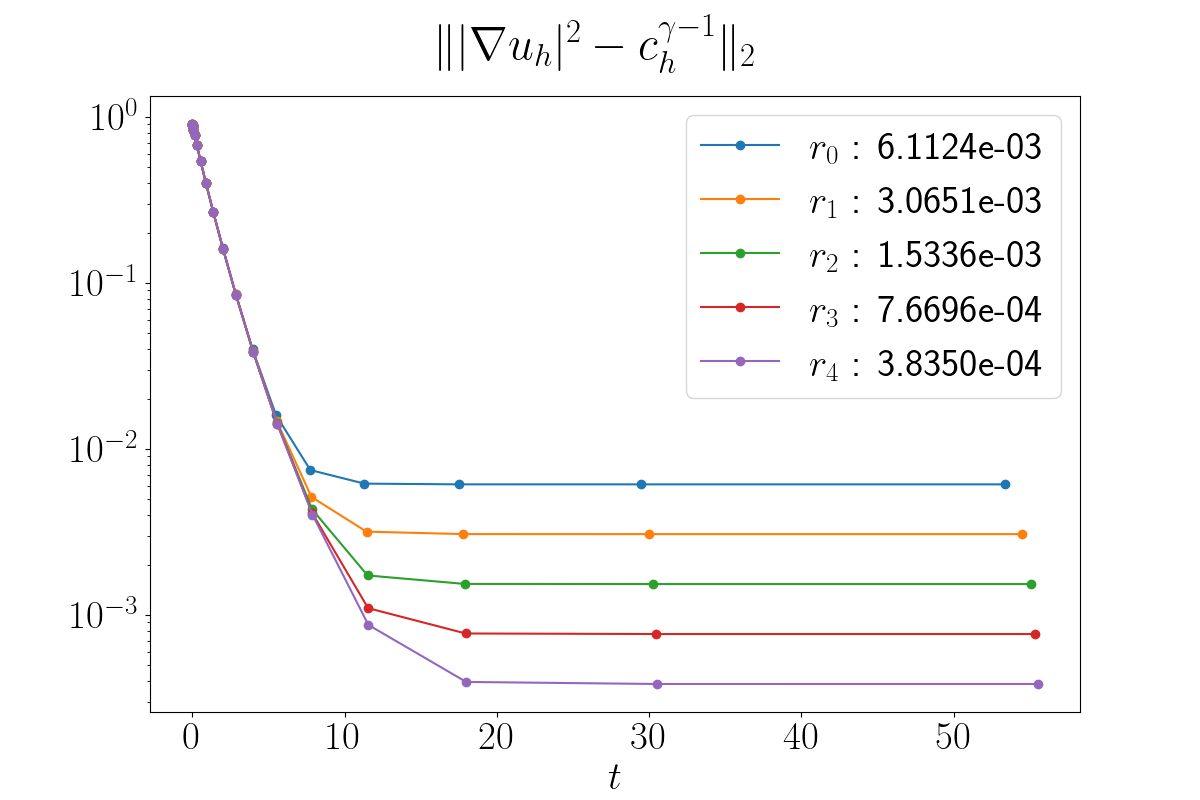} &
\includegraphics[width=0.22\textwidth,trim={50 0 70 0},clip]{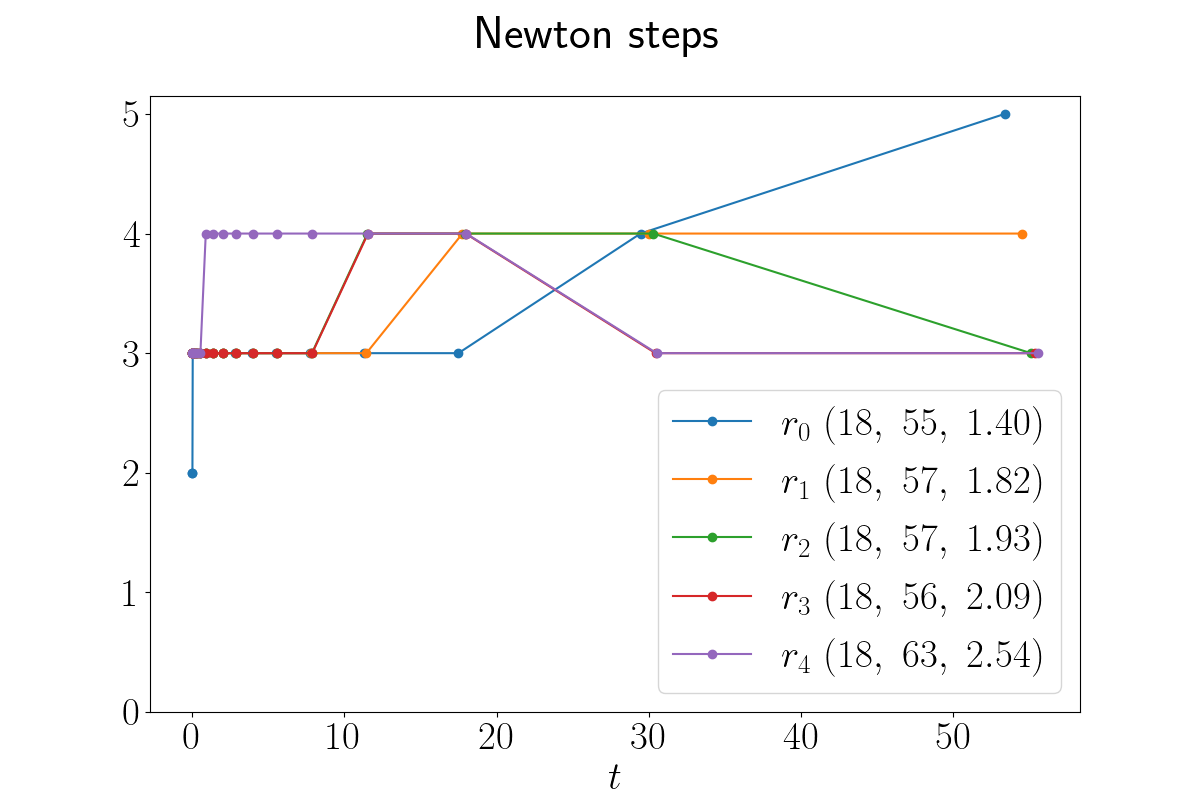} \\
& \multicolumn{2}{c}{TC4 $(p=20, \,\gamma=10/9)$} & \\
\includegraphics[width=0.22\textwidth,trim={50 0 70 0},clip]{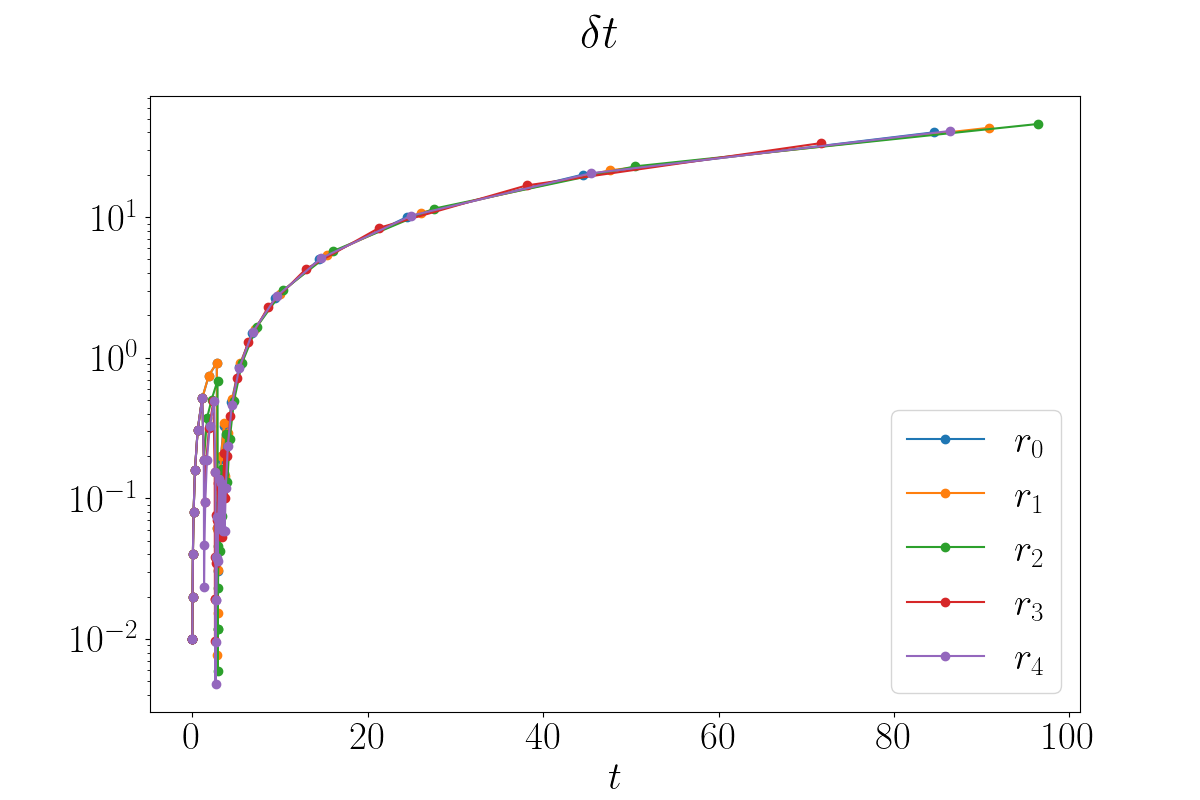} &
\includegraphics[width=0.22\textwidth,trim={50 0 70 0},clip]{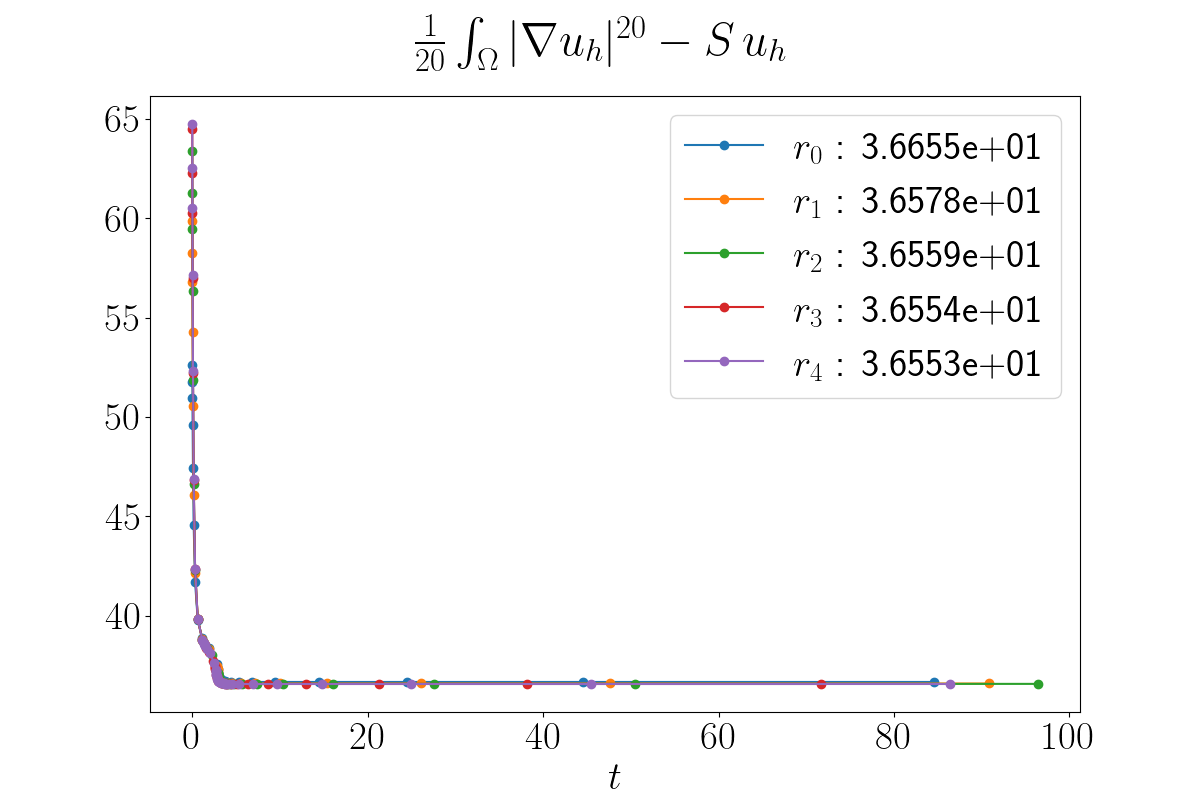} &
\includegraphics[width=0.22\textwidth,trim={50 0 70 0},clip]{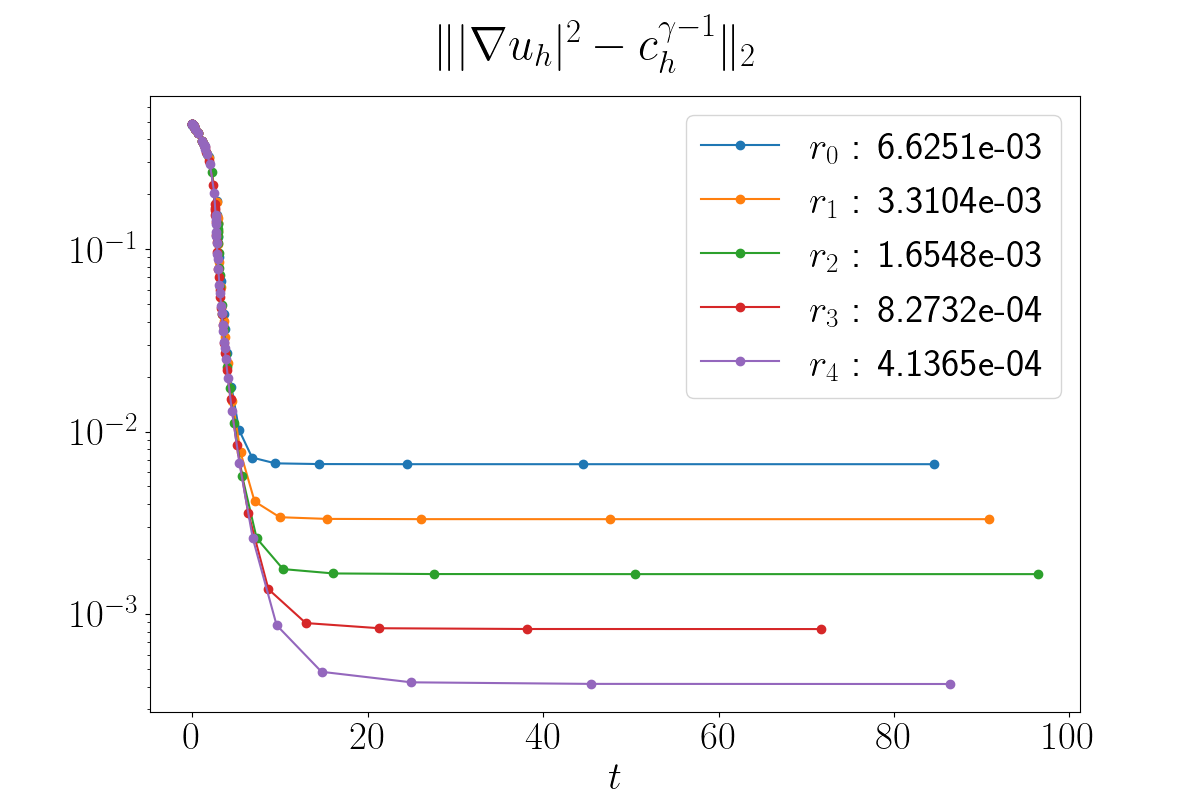} &
\includegraphics[width=0.22\textwidth,trim={50 0 70 0},clip]{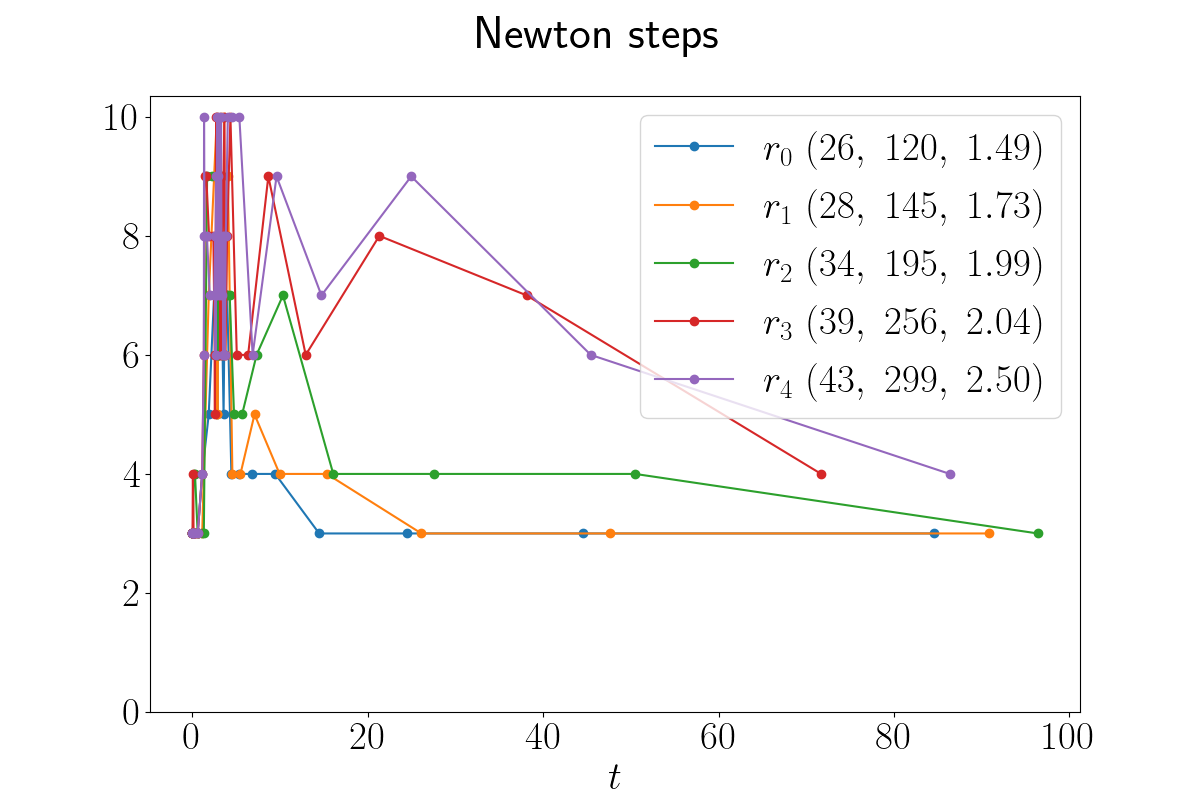} \\
\end{tabular}
\caption{From left to right: time step, $p$-Laplacian energy $1/p \int_\Omega |\nabla u_h|^p - S u_h$ (final value in the legend), $L^2$ norm of $|\nabla u_n|^2 - c_n^{\gamma-1}$, and number of Newton steps as a function of time for different levels of uniform refinement $r$ for TC2 (top row) and TC4 (bottom row). The legends in the right panels report in parentheses the total number of time steps, the total number of nonlinear steps, and the average number of Krylov iterations per Newton step.}
\label{fig:plap_bl3_solver}
\end{figure}

Details on the solution algorithm for the cases TC2 and TC4 are reported in Figure \ref{fig:plap_bl3_solver}; for each time step $n$, we report, from left to right, the time increment $\delta t_n$, the $p$-Laplacian energy $\frac{1}{p} \int_\Omega |\nabla u_h|^p - S u_h$, the $L^2$ norm of $|\nabla u_n|^2 - c_n^{\gamma-1}$, and the number of Newton steps as a function of simulation time for the different levels of refinement, from the coarsest $r_0$, to the finest $r_4$. The legends in the right panels collect the total number of time steps, the total number of Newton steps, and the average number of Krylov iterations per Newton step to provide additional information on the computational costs of the method.

In both the test cases, the time step is adaptively adjusted to follow the relaxation dynamics, which is smoother in the TC2 case. The $p$-Laplacian energy and $\| |\nabla u_n|^2 - c_n^{\gamma-1}\|_2$ decrease over time, with $\| |\nabla u_n|^2 - c_n^{\gamma-1}\|_2$ being reduced up to the discretization error $\bigO(h)$ because we evaluate the quantity at quadrature points and $\nabla u_n \notin \Cspace_h$ with quadrilateral elements. 
For TC2, the number of time steps and the total number of Newton steps remain independent of the refinement level, making the convergence of the method mesh-independent. The linear preconditioning strategy is also highly effective, requiring, on average, only two Krylov iterations per Newton step in all cases. For the TC4 case, time relaxation is more challenging, and a larger value of $p$ leads to stiffer dynamics at all refinement levels. After an initial transient phase, during which the Poisson problem with constant coefficients relaxes smoothly, the increasing nonlinearity of the problem causes the method to enter a second phase that demands a higher number of both time and Newton steps. This occurs because we cannot guarantee the non-negativity of $c_h$ within the nonlinear time step solver, which in turn fails more frequently, forcing the time step to shrink. Nevertheless, as the mesh is refined, these numbers appear to plateau. To address these challenges, we plan to study modifications of our method considering a mixed formulation of the Poisson problem, using an exponential map to handle the conductivity, and investigating alternative time-advancing approaches such as pseudo-transient continuation \cite{Coffey2003}.

\begin{figure}[htbp]
\centering
\includegraphics[width=0.3\textwidth,trim={50 0 70 50},clip]{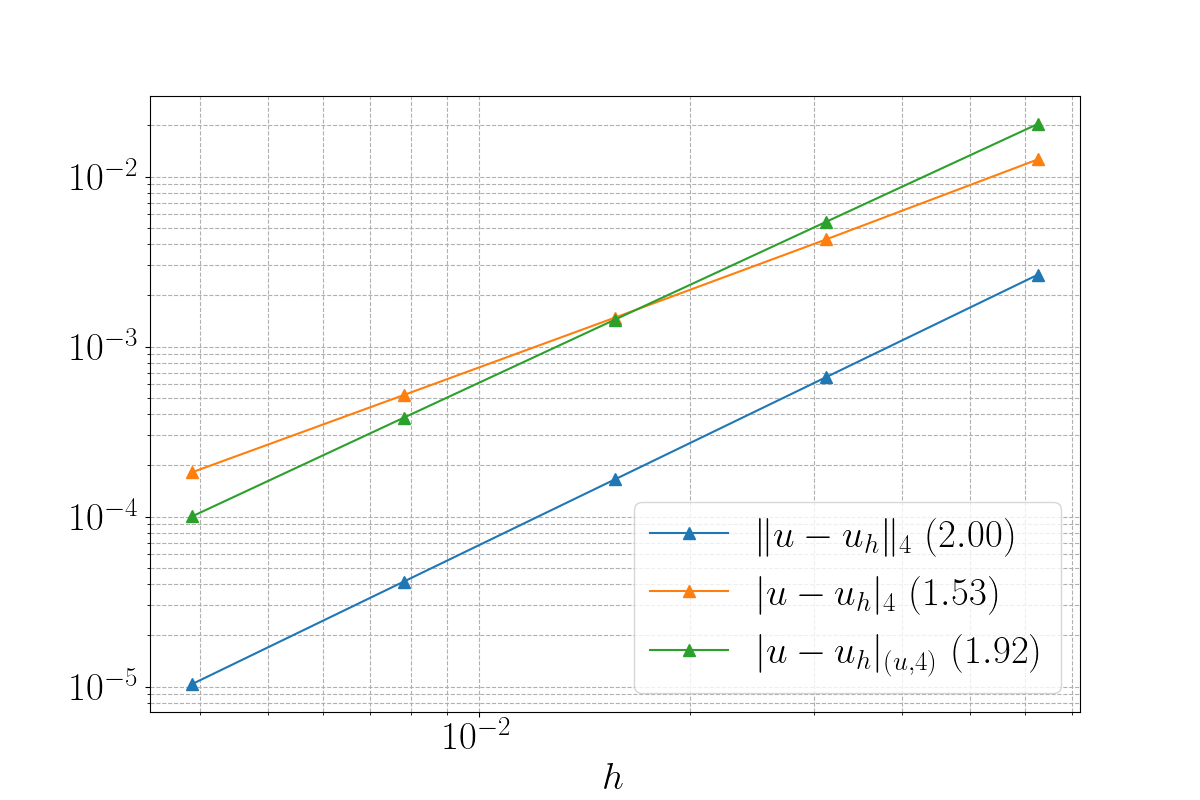}
\caption{TC5: MMS errors and convergence rates (in parentheses) for different error metrics.}
\label{fig:plap_bl5_conv}
\end{figure}

We now consider a more challenging test case with a non-convex domain and mixed boundary conditions for the potential, using an analytical solution expressed in polar coordinates
\begin{equation}\label{eq:ckeasy}
\begin{split}
u(r,\theta) &= r^\alpha \sin(\alpha \theta),\\
S(r,\theta) &= -\alpha  |\alpha| ^ {p - 2} (\alpha - 1)(p-2)r^{p\alpha - \alpha - p}\sin(\alpha \theta),\\
\end{split}
\end{equation}
with $\alpha=2$ and $p=4$.
This test case, referred to as TC5, is defined on an L-shaped domain $\Omega=(-1,1)^2\setminus [0,1)\times(-1,0]$. We impose non-homogeneous essential boundary conditions for the potential on $\Gamma_D=(\{0\}\times[-1,0])\cup([0,1]\times\{0\})$ where the prescribed values are taken from the analytical solution. On the complementary part of the boundary $\Gamma_N = \Gamma \setminus \Gamma_D$ we enforce non-homogeneous natural boundary conditions by adding a boundary integral to the right-hand side of the Poisson equation
\begin{equation*}
\int_\Omega(c_n + r)\nabla u_n \cdot \nabla \ptest_i ~ \dx =  \int_\Omega S \ptest_i ~ \dx +\int_{\Gamma_N} |\nabla u|^{p-2} \frac{\partial u}{\partial n}\ptest_i ~ \ds,
\end{equation*}
where $u$ is the exact solution. Convergence rates with different levels of uniform refinement are given in Figure \ref{fig:plap_bl5_conv}, showing second-order convergence in the $L^p$ norm and superlinear convergence in the $W^{1,p}(\Omega)$ semi-norm and in the quasi-norm $|\cdot|_{(u,p)}$.

\begin{figure}[htbp]
\centering
\begin{tabular}{c c c}
\includegraphics[width=0.3\textwidth,trim={50 0 70 0},clip]{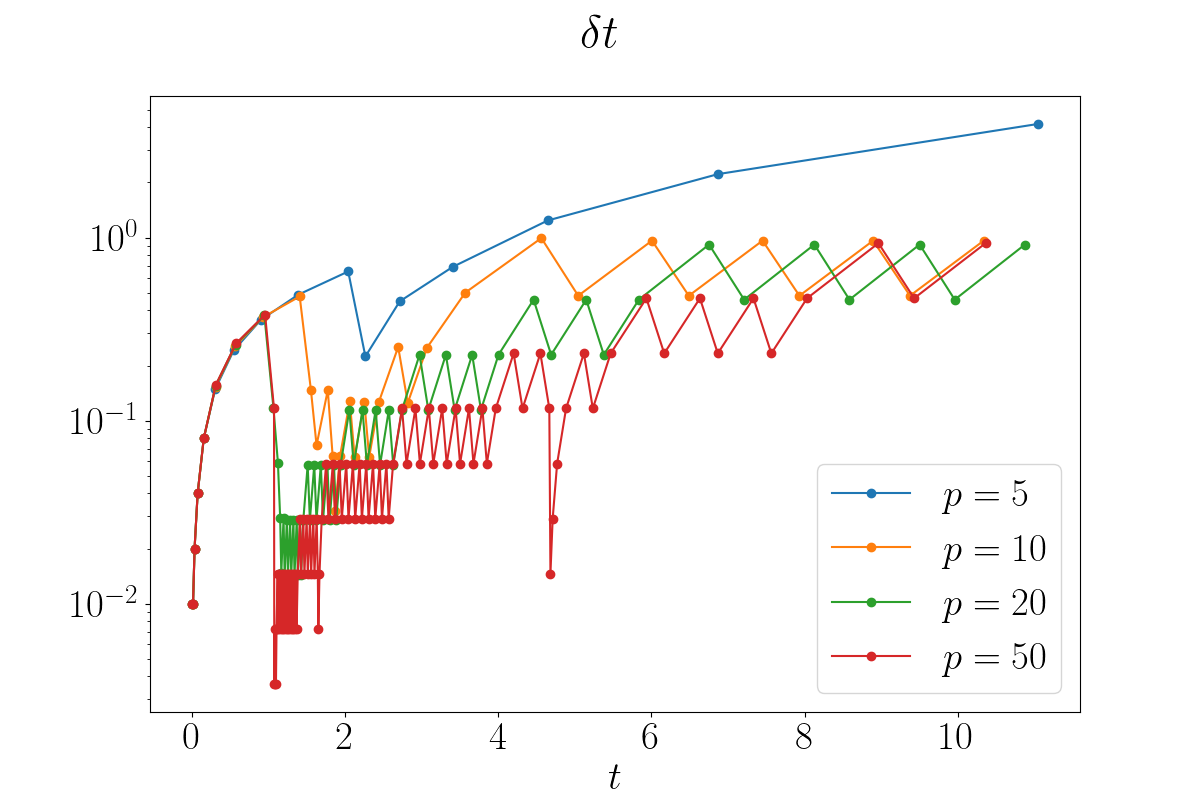} &
\includegraphics[width=0.3\textwidth,trim={50 0 70 0},clip]{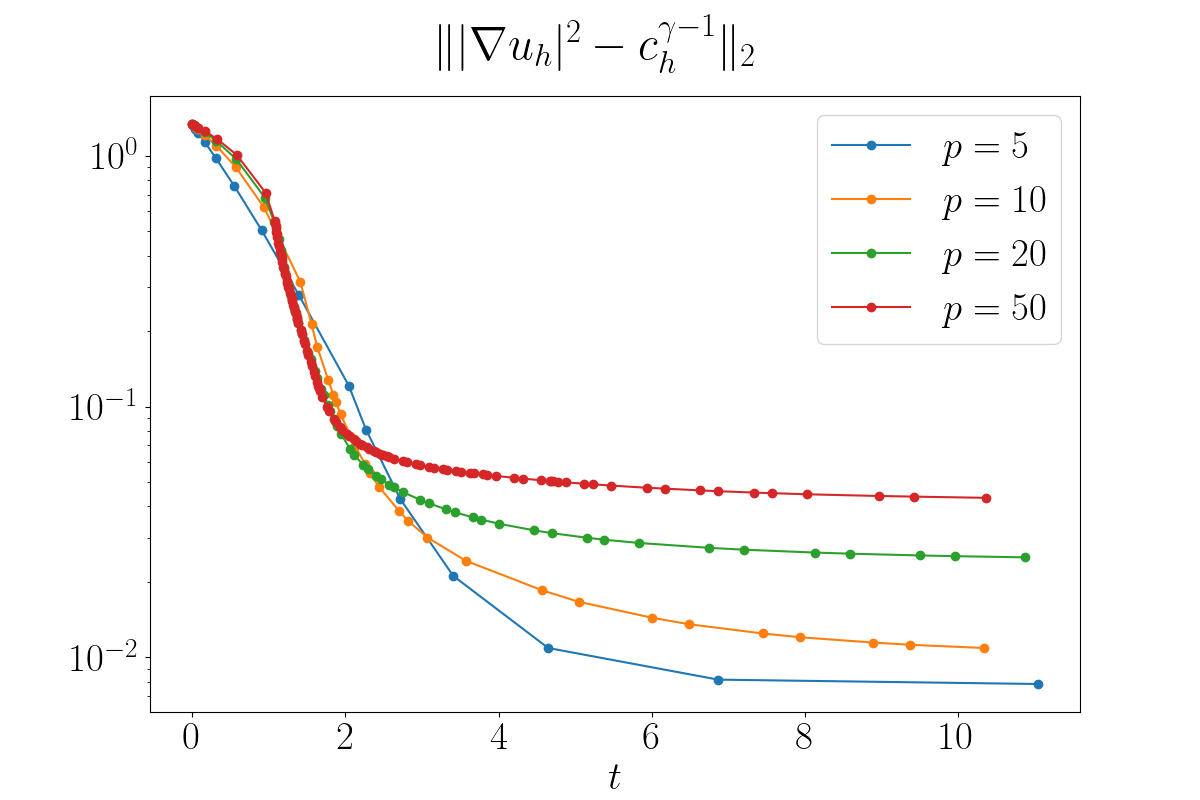} &
\includegraphics[width=0.3\textwidth,trim={50 0 70 0},clip]{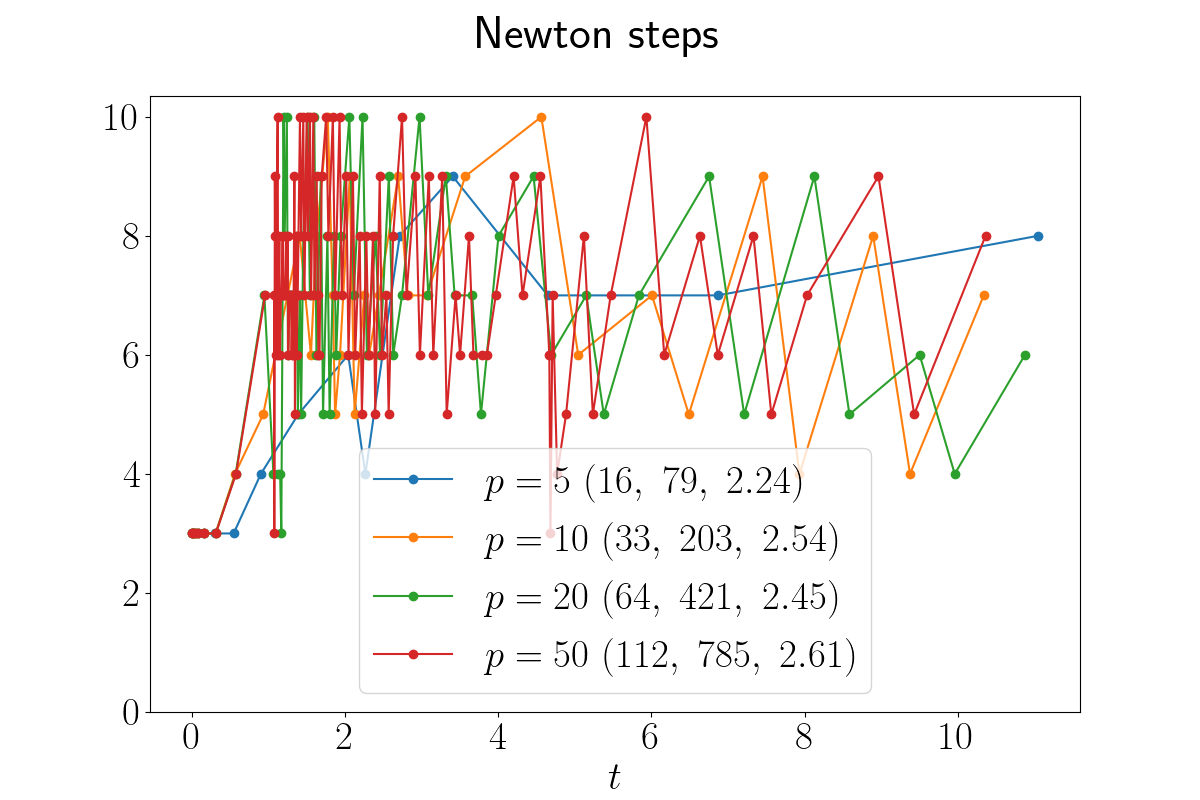} \\
\end{tabular}
\caption{TC6: Time step (left panel), $L^2$ norm of $|\nabla u_n|^2 - c_n^{\gamma-1}$ (center), and number of Newton steps (right) as a function of time for different $p$-Laplacian exponents $p$. The legends in the right panels report in parentheses the total number of time steps, the total number of nonlinear steps, and the average number of Krylov iterations per Newton step.}
\label{fig:plap_lshaped_solver}
\end{figure}

We conclude this Section by presenting results for different values of $p$ in a challenging problem without an analytical solution, following the experimental setting of Section 7.3 in \cite{BalciDiening2023}. This test, referred to as TC6, considers a constant source term $S=2$ and homogeneous essential boundary conditions for the potential on the L-shaped domain $(-1,1)^2 \setminus [0,1)\times (-1,0]$. The domain is uniformly discretized using ten thousand elements, with mesh size $h = 1.9 \times 10^{-3}$. Solver performance for $p$-Laplacian exponents $p=5,\,10,\,20,\,50$ (corresponding to $\gamma=5/3,\,5/4,\,10/9,\,25/24$) is reported in Figure \ref{fig:plap_lshaped_solver}.

\begin{figure}[htbp]
\centering
\begin{tabular}{c c}
$u_h$ & $c_h$\\
\includegraphics[width=0.4\textwidth,trim={70 30 80 50},clip]{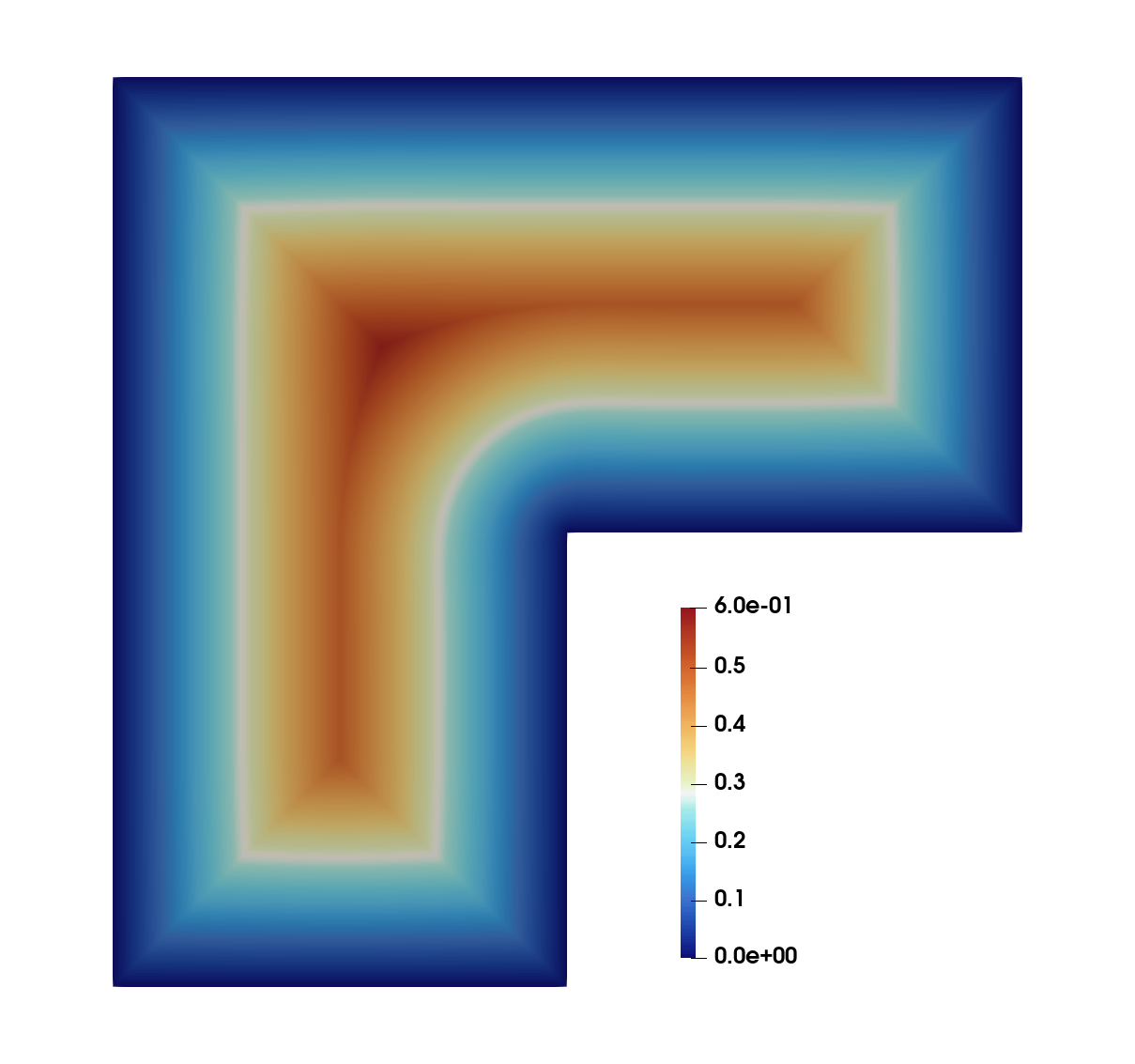} &
\includegraphics[width=0.4\textwidth,trim={70 30 80 50},clip]{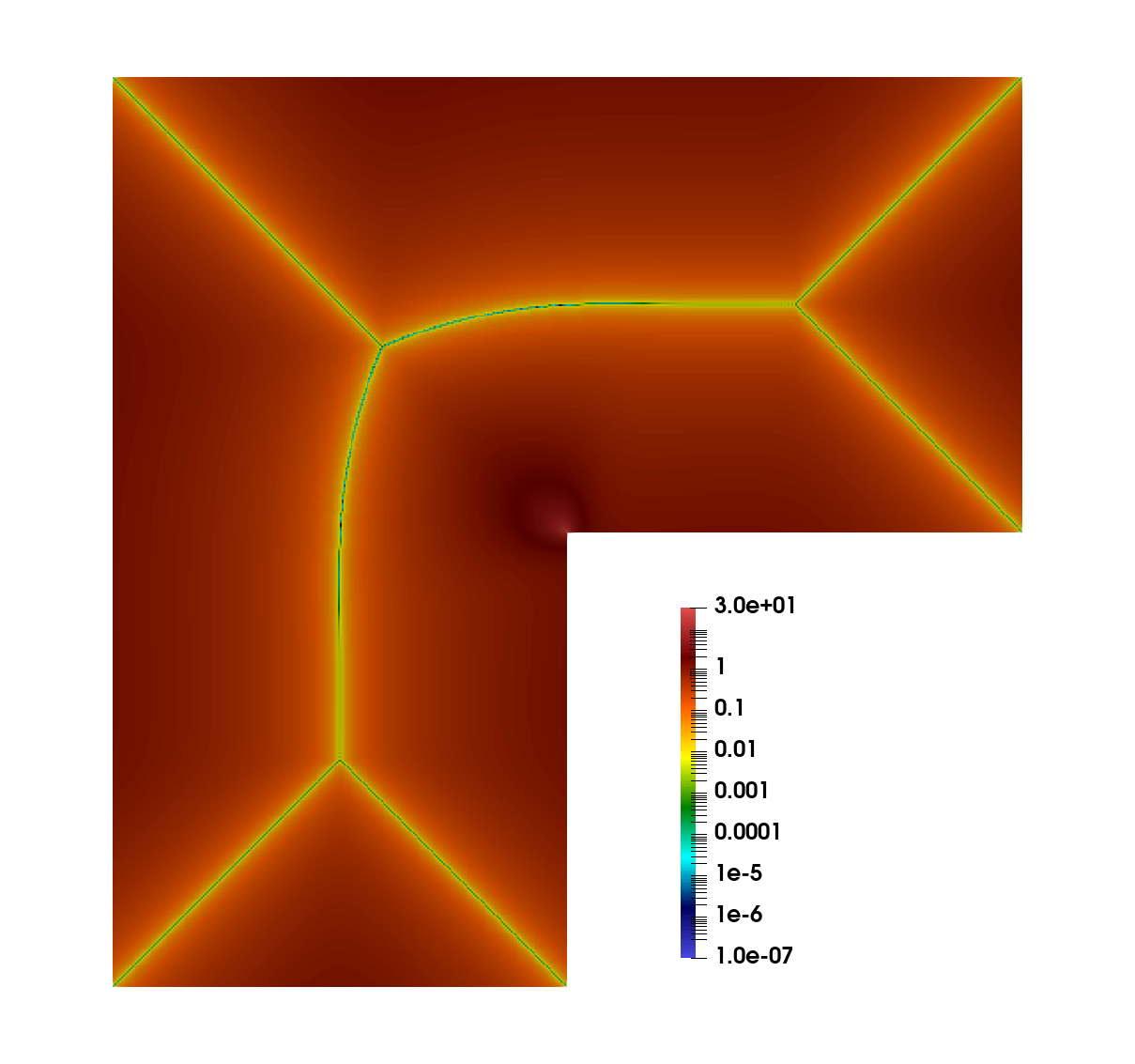}\\ 
$|\nabla u_h|^2 - c_h^{\gamma-1}$ & $\pi_h(|\nabla u_h|^2 - c_h^{\gamma-1})$\\
\includegraphics[width=0.4\textwidth,trim={70 30 80 50},clip]{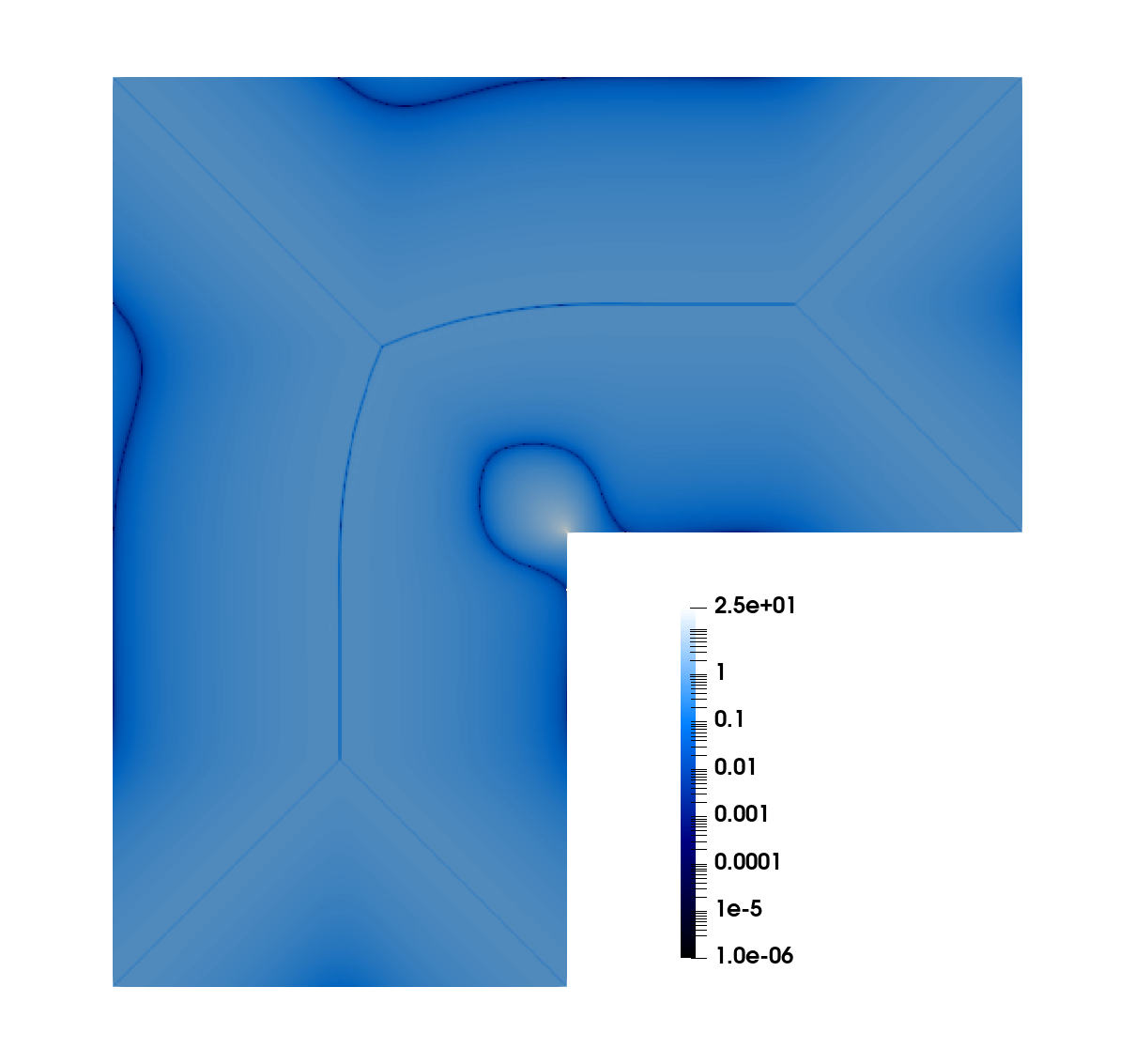} &
\includegraphics[width=0.4\textwidth,trim={70 30 80 50},clip]{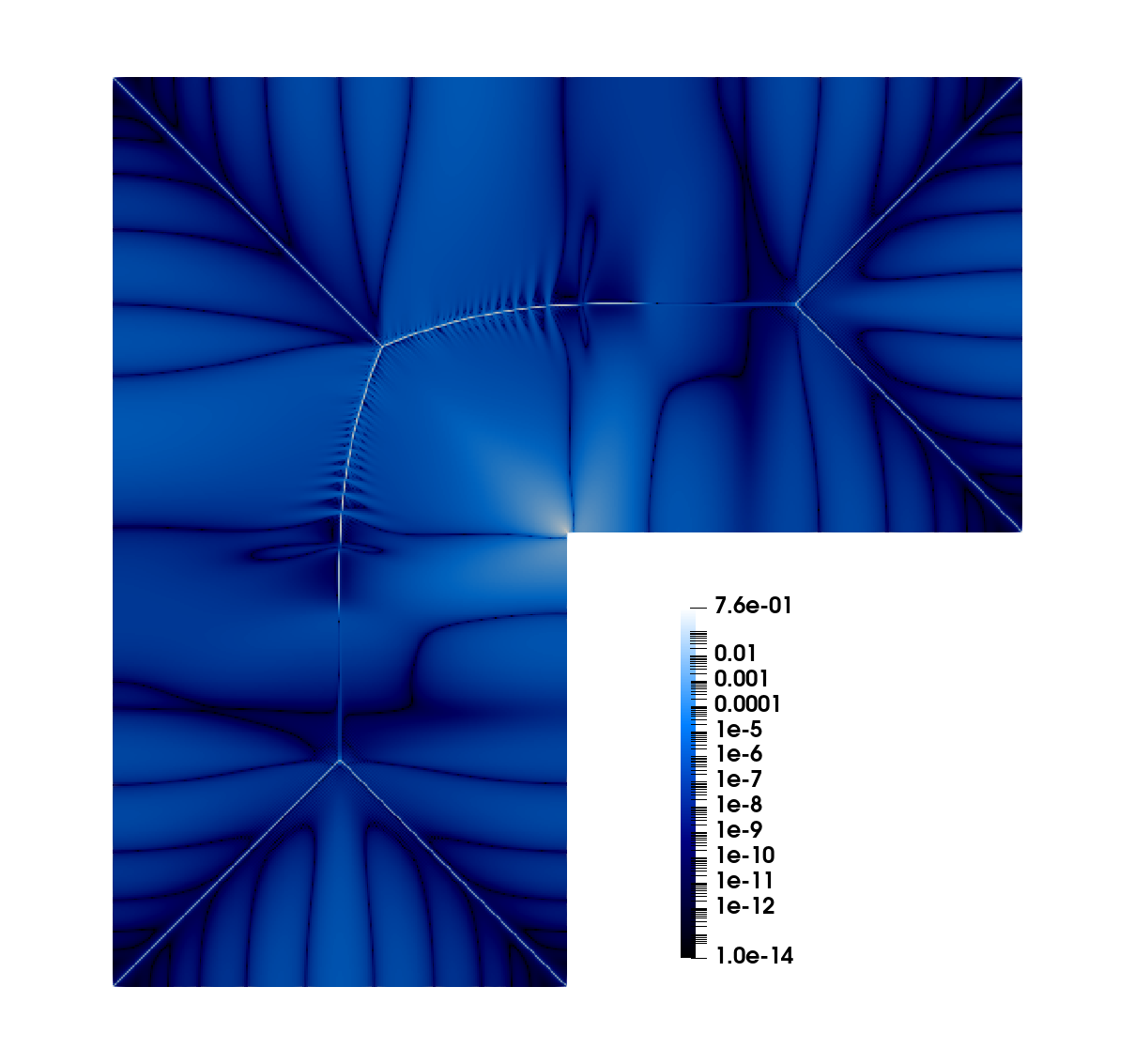} 
\end{tabular}
\caption{Top row: final state for potential (left) and conductivity (right) for the TC6 test case with $p=50$ (corresponding to $\gamma=25/24$). Bottom row: the post-processed final point-wise error for $|\nabla u_h|^2 - c_h^{\gamma-1}$ (left) and its projection onto $L^2(\Omega)$ (right).}
\label{fig:plap_lshaped_p50}
\end{figure}

The final values of $\| |\nabla u_n|^2 -  c_n^{\gamma-1}\|_2$ are larger for larger $p$, possibly because we did not use a $p$-dependent norm. The solver is not $p$-independent overall, since it requires more time and Newton steps as $p$ is increased, as a result of the increase in nonlinear solver failures. However, the average number of Krylov iterations per Newton step remains essentially constant across all cases, confirming the robustness of the Schur complement-based preconditioner. 
Figure \ref{fig:plap_lshaped_p50} shows the final states of the potential $u_h$ and the conductivity $c_h$ for $p=50$, together with the pointwise error $|\nabla u_h|^2 - c_h^{\gamma-1}$ (evaluated in post-processing) and its projection onto $L^2(\Omega)$, $\pi_h\big(|\nabla u_h|^2 - c_h^{\gamma-1}\big)$ (evaluated within the finite element simulation code). The conductivity (plotted on a logarithmic scale) is small in regions of sharp potential gradients, where the singularities appear to have Hausdorff dimension one, and attains large values near the concave corner of the domain. These patterns closely resemble the mesh adaptation patterns reported in Figure 7.3 of \cite{BalciDiening2023}. The error distribution is nearly uniform across the domain and noticeably larger near the L-shaped corner; the $L^2$ projection of the error is small everywhere in the domain and slightly larger near the corner and where $c_h \approx 0$. Future work will focus on adaptive mesh refinement strategies and the development of error indicators specifically tailored to the differential-algebraic transient dynamics of the solver.

\end{document}